\documentclass[final]{siamltex}
\usepackage{graphicx,enumerate,amsmath,amssymb,afterpage,cite}

\def\figdir{.}

\newcommand{\dif}{\,\mathrm{d}}

\newcommand{\parderiv}[2]{\frac{\partial #1}{\partial #2}}
\newcommand{\deriv}[2]{\frac{\mathrm{d} #1}{\mathrm{d} #2}}
\newcommand\BS\boldsymbol
\newcommand\picturehere[1]{\includegraphics[width=0.5\textwidth]{#1}}

\newcommand\comment[1]{}

\title{A simple closure approximation for slow dynamics of a
  multiscale system: nonlinear and multiplicative coupling}

\author{Rafail V. Abramov\thanks{Department of Mathematics,
Statistics and Computer Science, University of Illinois at Chicago
({\tt abramov@math.uic.edu})}}

\begin{document}

\maketitle

\begin{abstract}
Multiscale dynamics are ubiquitous in applications of modern
science. Because of time scale separation between relatively small set
of slowly evolving variables and (typically) much larger set of
rapidly changing variables, direct numerical simulations of such
systems often require relatively small time discretization step to
resolve fast dynamics, which, in turn, increases computational
expense. As a result, it became a popular approach in applications to
develop a closed approximate model for slow variables alone, which
both effectively reduces the dimension of the phase space of dynamics,
as well as allows for a longer time discretization step. In this work
we develop a new method for approximate reduced model, based on the
linear fluctuation-dissipation theorem applied to statistical states
of the fast variables. The method is suitable for situations with
quadratically nonlinear and multiplicative coupling. We show that,
with complex quadratically nonlinear and multiplicative coupling in
both slow and fast variables, this method produces comparable
statistics to what is exhibited by an original multiscale model. In
contrast, it is observed that the results from the simplified closed
model with a constant coupling term parameterization are consistently
less precise.
\end{abstract}



\pagestyle{myheadings}
\thispagestyle{plain}
\markboth{Rafail V. Abramov}{A closure for slow
dynamics of a multiscale system: nonlinear coupling}

\section{Introduction}

Multiscale dynamics are ubiquitous in applications of modern science,
with geophysical science and climate change prediction being
well-known examples \cite{FraMajVan,Has,BuiMilPal,Pal3}. Direct
numerical simulations of multiscale systems are difficult, both
because a relatively short time discretization step is required to
resolve fast dynamics, and due to a large number of dynamical
variables. Moreover, in some applications such as climate change
prediction one is interested in long-term statistics of slow dynamics,
which further increases computational expense.

A popular approach for simulating multiscale dynamics in practice with
limited computational resources is to create an approximate reduced
model for slow variables alone, which allows to increase the length of
the time discretization step and reduced the dimension of the phase
space of the system, which is accomplished via an approximate closure
of the coupling terms between slow and fast variables of the system.
Many closure methods were designed for multiscale dynamical systems
\cite{CroVan,FatVan,MajTimVan,MajTimVan2,MajTimVan3,MajTimVan4}, based
on the averaging formalism for the fast variables
\cite{Pap,Van,Vol}. Some methods approximate the coupling terms with
appropriate stochastic processes
\cite{MajTimVan,MajTimVan2,MajTimVan3,MajTimVan4,Wilks} or conditional
Markov chains \cite{CroVan}, while others \cite{FatVan} parameterize
slow-fast interactions by direct tabulation and curve fitting.

In a recent work \cite{Abr9} the author developed a simple approach of
computing the reduced model for slow variables alone via a single
computation of relevant statistics for the fast dynamics with a fixed
reference state of the slow variables, using the linear
fluctuation-dissipation theorem (FDT)
\cite{Abr5,Abr6,Abr7,Abr8,AbrMaj4,AbrMaj5,AbrMaj6,AbrMaj7,MajAbrGro,Ris}.
It was shown that, for the appropriately rescaled two-scale Lorenz 96
model \cite{Abr8} with linear coupling between slow and fast
variables, the method reproduced statistics of the slow variables of
the complete two-scale Lorenz model with good precision. However,
nonlinear coupling was not addressed in \cite{Abr9}.

This work is a natural extension of the method in \cite{Abr9} onto
quadratically nonlinear and multiplicative types of coupling in both
slow and fast variables. It is also based on the linear FDT, which is
used to compute the response of linear, nonlinear and multiplicative
terms of the fast variables coupling to changes in the slow
variables. We show numerical experiments for this method with the
two-scale Lorenz model and general forms of coupling which include
multiplicative and quadratically nonlinear terms.

The manuscript is organized as follows. In Section \ref{sec:theory} we
formulate the theory for the new method. In Section \ref{sec:lorenz}
we present the two-scale Lorenz model \cite{FatVan,Lor,LorEma},
rescaled in such a way that the mean states and variances of both the
fast and slow variables are near zero and one, respectively
\cite{Abr8,Abr9}. Section \ref{sec:num} shows the results of numerical
simulations with both the two-scale Lorenz model and the reduced model
for slow variables only, comparing different statistics of the time
series. Section \ref{sec:sum} summarizes the results of this work.

\section{Derivation of the reduced model}
\label{sec:theory}

We start with a general two-scale system of differential equations of
the form
\begin{equation}
\label{eq:dyn_sys}
\deriv{\BS x}t=\BS F(\BS x,\BS y),\qquad \deriv{\BS y}t=\BS G(\BS
x,\BS y),
\end{equation}
where $\BS x=\BS x(t)\in\mathbb R^{N_x}$ are the slow variables, $\BS
y=\BS y(t)\in\mathbb R^{N_y}$ are the fast variables, and $\BS F$ and
$\BS G$ are $N_x$ and $N_y$ vector-valued functions of $\BS x$ and
$\BS y$, respectively. We assume that the $\BS y$-variables are
sufficiently fast for a valid approximation of the system in
\eqref{eq:dyn_sys} by the averaged dynamics for $\BS x$, given by
\begin{equation}
\label{eq:dyn_sys_slow_limiting_x}
\deriv{\BS x}t=\langle\BS F\rangle(\BS x),\qquad \langle\BS
F\rangle(\BS x)=\int_{\mathbb R^{N_y}}\BS F(\BS x,\BS z) \dif\mu_{\BS
  x}(\BS z),
\end{equation}
for finite times (for a more detailed description of the averaging
formalism, see \cite{Abr6,Abr8,Pap,Van,Vol}). Here, $\mu_{\BS x}$
denotes the invariant probability measure of the limiting fast
dynamics, which are given by system
\begin{equation}
\label{eq:dyn_sys_fast_limiting_z}
\deriv{\BS z}\tau=\BS G(\BS x,\BS z).
\end{equation}
Above, $\BS x$ is a constant parameter, and the solution of
\eqref{eq:dyn_sys_fast_limiting_z} is given by the flow $\BS
z(\tau)=\phi_{\BS x}^\tau\BS z_0$. We tacitly assume that all typical
initial conditions $\BS z_0$ fall into the support of the same ergodic
component of $\mu_{\BS x}$, and that $\langle\BS F\rangle(\BS x)$
varies smoothly with respect to $\BS x$, as it often happens when
$\mu_{\BS x}$ is an SRB measure \cite{EckRue,Rue1,Rue2,Rue3,You}.
Using the ergodicity assumption for $\mu_{\BS x}$, we can practically
compute the measure average via the time average
\begin{equation}
\label{eq:F_time_average}
\langle\BS F\rangle(\BS x)=\lim_{r\to\infty}\frac 1r\int_0^r\BS F(\BS
x, \BS z(\tau))\dif\tau,
\end{equation}
where $\BS z(\tau)$ is a long-term trajectory of
\eqref{eq:dyn_sys_fast_limiting_z}.

Following \cite{Abr9}, here we propose an approximation to $\langle\BS
F\rangle(\BS x)$, based on the linear FDT. In order to do this,
certain assumptions have to be made regarding coupling, that is the
dependence of $\BS F$ and $\BS G$ on the fast and slow variables,
respectively. In \cite{Abr9} it was assumed that $\BS F$ depends
linearly on the fast variables, and $\BS G$ depends linearly on the
slow variables (linear coupling). Here we assume that $\BS F$ is
quadratic in the fast variables in the vicinity of $\BS{\bar z}(\BS
x)$, which is the mean state of \eqref{eq:dyn_sys_fast_limiting_z}
with $\BS x$ set as a constant parameter:
\begin{equation}
\begin{split}
\BS F(\BS x,\BS z)&=\BS F(\BS x,\BS{\bar z}(\BS x))+ \parderiv{\BS
  F}{\BS y}(\BS x,\BS{\bar z}(\BS x))(\BS z-\BS{\bar z}(\BS x))+\\&+
\frac 12\parderiv{^2\BS F}{\BS y^2}(\BS x,\BS{\bar z}(\BS x))
:(\BS z-\BS{\bar z}(\BS x))\otimes(\BS z-\BS{\bar z}(\BS x)),
\end{split}
\end{equation}
where $\partial\BS F/\partial\BS y$ and $\partial^2\BS F/\partial\BS
y^2$ denote the first and second partial derivatives of $\BS F$ with
respect to its second argument, respectively, and ``:'' is the
element-wise (Hadamard) matrix product with summation. The assumption
of quadratic dependence of $\BS F$ on the fast variables is not too
restrictive for practical applications; indeed, many real-world
geophysical processes have at most quadratic dependence on velocity or
streamfunction fields due to advection. Now, the average with respect
to $\mu_{\BS x}$ is given by
\begin{equation}
\label{eq:averaged_F}
\langle\BS F\rangle(\BS x)=\BS F(\BS x,\BS{\bar z}(\BS x))+\frac 12
\parderiv{^2\BS F}{\BS y^2}(\BS x,\BS{\bar z}(\BS x)):\BS\Sigma(\BS x),
\end{equation}
where $\BS\Sigma(\BS x)$ is the covariance of $\BS z$, centered at
$\BS{\bar z}(\BS x)$. The mean state $\BS{\bar z}(\BS x)$ and the
covariance $\BS\Sigma(\BS x)$ are given by
\begin{subequations}
\begin{equation}
\BS{\bar z}(\BS x)=\int_{\mathbb R^{N_y}}\BS z\dif\mu_{\BS x}(\BS z),
\end{equation}
\begin{equation}
\BS\Sigma(\BS x)=\int_{\mathbb R^{N_y}}(\BS z-\BS{\bar z}(\BS
x))\otimes(\BS z-\BS{\bar z}(\BS x))\dif\mu_{\BS x}(\BS z).
\end{equation}
\end{subequations}
As we can see, the average of $\BS F(\BS x,\BS z)$ with respect to
$\mu_{\BS x}$ is now expressed as a nonlinear function in the mean
state $\BS{\bar z}(\BS x)$, and linear function in covariance
$\BS\Sigma(\BS x)$. If we know how these quantities respond to changes
in $\BS x$, we can also calculate the approximation of the $\BS
x$-dependent average $\langle\BS F\rangle(\BS x)$.

Now, the response of $\BS{\bar z}(\BS x)$ and $\BS\Sigma(\BS x)$ to
changes in $\BS x$ can be estimated via the linear FDT. In order to
obtain the response formulas, we have to impose some restrictions on
the structure of $\BS G(\BS x,\BS y)$ in $\BS y$. Here, we assume that
$\BS G(\BS x,\BS y)$ can be written as
\begin{equation}
\BS G(\BS x,\BS y)=\BS g(\BS y)+\BS H(\BS x)\BS y+\BS h(\BS x),
\end{equation}
where $\BS g(\BS y)$ is a $N_y$-vector nonlinear function of $\BS y$,
$\BS h(\BS x)$ is a $N_y$-vector function of $\BS x$, and $\BS H(\BS
y)$ is a $N_y\times N_y$-matrix valued function of $\BS x$. Again,
this assumption is not too restrictive for practical applications, as
only the nonlinear part of $\BS G$ does not depend on $\BS x$. In this
case, the limiting system in \eqref{eq:dyn_sys_fast_limiting_z} can be
written as
\begin{equation}
\deriv{\BS z}{\BS\tau}=\BS g(\BS z)+\BS H(\BS x)\BS z+\BS h(\BS
x),
\end{equation}
where $\BS x$ is a constant parameter. Now, let $\BS H^*=\langle\BS
H(\BS x)\rangle$, $\BS h^*=\langle\BS h(\BS x)\rangle$ denote the
long-term averages of $\BS H(\BS x)$ and $\BS h(\BS x)$ over a
trajectory of \eqref{eq:dyn_sys}, so that we can write the fast
limiting system as
\begin{equation}
\label{eq:dyn_sys_fast_limiting_delta}
\begin{split}
\deriv{\BS z}{\BS\tau}=\BS g(\BS z)+(\BS H^*+\delta\BS H(\BS x))
\BS z+(\BS h^*+\delta\BS h(\BS x)),\\
\delta\BS H(\BS x)=\BS H(\BS x)-\BS H^*,\qquad\delta\BS h(\BS
x)=\BS h(\BS x)-\BS h^*.
\end{split}
\end{equation}
Given the mean state $\BS{\bar z}^*$ and the mean-centered covariance
matrix $\BS\Sigma^*$ for unperturbed
\eqref{eq:dyn_sys_fast_limiting_delta} with $\delta\BS H(\BS x)$ and
$\delta\BS h(\BS x)$ set to zeros, one can think of $\BS{\bar z}(\BS
x)$ and $\BS\Sigma(\BS x)$ as responses of $\BS{\bar z}^*$ and
$\BS\Sigma^*$ to nonzero $\delta\BS H(\BS x)$ and $\delta\BS h(\BS
x)$. These responses can be written as linear approximations via the
FDT:
\begin{equation}
\label{eq:z_Sigma}
\begin{split}
\BS{\bar z}(\BS x)&\approx\BS{\bar z}^*+\BS R^{h\to\bar z}(\delta\BS h(\BS x)
+\delta\BS H(\BS x)\BS{\bar z}^*)+
\BS R^{H\to\bar z}\delta\BS H(\BS x),\\
\BS\Sigma(\BS x)&\approx\BS\Sigma^*+\BS R^{h\to\Sigma}(\delta\BS h(\BS x)
+\delta\BS H(\BS x)\BS{\bar z}^*)+
\BS R^{H\to\Sigma}\delta\BS H(\BS x),
\end{split}
\end{equation}
where the linear response operators $\BS R^{h\to\bar z}$, $\BS
R^{H\to\bar z}$, $\BS R^{h\to\Sigma}$, $\BS R^{H\to\Sigma}$ are given
by the quasi-Gaussian formulas
\begin{subequations}
\label{eq:R_Gaussian}
\begin{equation}
\BS R_{ij}^{h\to\bar z}=\int_0^\infty\left[\lim_{r\to\infty}\frac 1r
\int_0^r(z_i(t+s)-\bar z^*_i)(z_k(t)-\bar z^*_k)\dif t\right]\dif s
\;\Sigma^{*-1}_{kj},
\end{equation}
\begin{equation}
\begin{split}
\BS R_{ijk}^{H\to\bar z}=\int_0^\infty\bigg[\lim_{r\to\infty}\frac 1r
\int_0^r(z_i(t+s)-\bar z^*_i)(z_l(t)-\bar z^*_l)\times\\\times
(z_k(t)-\bar z^*_k)\dif t\bigg]\dif s\;\Sigma^{*-1}_{lj},
\end{split}
\end{equation}
\begin{equation}
\begin{split}
\BS R_{ijk}^{h\to\Sigma}=\int_0^\infty\bigg[\lim_{r\to\infty}
\frac 1r\int_0^r(z_i(t+s)-\bar z^*_i)(z_j(t+s)-\bar z^*_j)\times\\
\times(z_l(t)-\bar z^*_l)\dif t\bigg]\dif s\;\Sigma^{*-1}_{lk},
\end{split}
\end{equation}
\begin{equation}
\begin{split}
\BS R_{ijkl}^{H\to\Sigma}=\int_0^\infty\bigg[\lim_{r\to\infty}
\frac 1r\int_0^r(z_i(t+s)-\bar z^*_i)(z_j(t+s)-\bar z^*_j)\times\\
\times(z_m(t)-\bar z^*_m)(z_l(t)-\bar z^*_l)\dif t\;
\Sigma^{*-1}_{mk}-\Sigma^*_{ij}\delta_{kl}\bigg]\dif s.
\end{split}
\end{equation}
\end{subequations}
Above, the time averaging is performed over a single long-term
trajectory of the unperturbed fast limiting system in
\eqref{eq:dyn_sys_fast_limiting_delta}, with $\delta\BS h$ and
$\delta\BS H$ set to zero. The details of derivation of
\eqref{eq:z_Sigma} and \eqref{eq:R_Gaussian} are given in Appendix
\ref{sec:app_reduced}, along with relevant references.

Combining \eqref{eq:z_Sigma} with \eqref{eq:averaged_F} yields the
approximate reduced system for slow variables alone. Indeed, observe
that if the averages $\BS{\bar z}^*$, $\BS\Sigma^*$, and the response
operators $\BS R^{h\to\bar z}$, $\BS R^{H\to\bar z}$, $\BS
R^{h\to\Sigma}$, and $\BS R^{H\to\Sigma}$ are computed, then
\eqref{eq:z_Sigma} and, therefore, \eqref{eq:averaged_F} are known
explicitly for given parameter $\BS x$.

\section{The Lorenz model with nonlinear coupling}
\label{sec:lorenz}

The rescaled Lorenz model with nonlinear coupling is given by
\begin{subequations}
\label{eq:lorenz_two_scale}
\begin{equation}
\begin{split}
\dot x_i=x_{i-1}(&x_{i+1}-x_{i-2})+\frac 1{\beta_x}(\bar x(x_{i+1}-x_{i-2})
-x_i)+\frac{F_x-\bar x}{\beta_x^2}-\\&-
\frac{\lambda_y}J\sum_{j=1}^J\left[(a+bx_i)y_{i,j}+(c+dx_i)(y_{i,j}^2-1)
\right],
\end{split}
\end{equation}
\begin{equation}
\begin{split}
\dot y_{i,j}=\frac 1\varepsilon\bigg[y_{i,j+1}
(y_{i,j-1}&-y_{i,j+2})+\frac 1{\beta_y}(\bar y(y_{i,j-1}-y_{i,j+2})-y_{i,j})
+\frac{F_y-\bar y}{\beta_y^2}\bigg]+\\&+\frac{\lambda_x}\varepsilon
\left[(a+cy_{i,j})x_i+(b+dy_{i,j})(x_i^2-1)\right],
\end{split}
\end{equation}
\end{subequations}
with $1\leq i\leq N_x$ and $1\leq j\leq J$, such that $N_y=N_xJ$. The
model has periodic boundary conditions: $x_{i+N_x}=x_i$,
$y_{i,j+J}=y_{i+1,j}$, and $y_{i+N_x,j}=y_{i,j}$. The parameters
$(\bar x,\beta_x)$ and $(\bar y,\beta_y)$ are the (mean,standard
deviation) pairs for the corresponding uncoupled and unrescaled Lorenz
models
\begin{subequations}
\begin{equation}
\dot x_i=x_{i-1}(x_{i+1}-x_{i-2})-x_i+F_x,
\end{equation}
\begin{equation}
\dot y_{i,j}=y_{i,j+1}(y_{i,j-1}-y_{i,j+2})-y_{i,j}+F_y,
\end{equation}
\end{subequations}
with the same periodic boundary conditions. The rescaling above
ensures that the Lorenz model in \eqref{eq:lorenz_two_scale} has zero
mean state and unit standard deviation for both slow and fast
variables in the absence of coupling ($\lambda_x=\lambda_y=0$), and
remain near these values when $\lambda_x$ and $\lambda_y$ are nonzero.
The original rescaled Lorenz model in \cite{Abr8,Abr9} with linear
coupling corresponds to the set of parameters $a=1$, $b=c=d=0$. The
nonlinear coupling above preserves the energy of the form
\begin{equation}
E=\frac{\lambda_x}{2\varepsilon}\sum_{i=1}^{N_x}x_i^2+\frac{\lambda_y}{2J}
\sum_{i=1}^{N_x}\sum_{j=1}^Jy_{i,j}^2.
\end{equation}
Indeed, observe that
\begin{equation}
\begin{split}
\deriv Et=&\frac{\lambda_x}\varepsilon\sum_{i=1}^{N_x}x_i\dot x_i+
\frac{\lambda_y}J\sum_{i=1}^{N_x}\sum_{j=1}^Jy_{i,j}\dot y_{i,j}=\\=
&-\frac{\lambda_x\lambda_y}{\varepsilon J}\sum_{i=1}^{N_x}\sum_{j=1}^J
\left(ax_iy_{i,j}+bx_i^2y_{i,j}+cx_iy_{i,j}^2+dx_i^2y_{i,j}^2\right)+\\&+
\frac{\lambda_x\lambda_y}{\varepsilon J}\sum_{i=1}^{N_x}\sum_{j=1}^J
\left(ax_iy_{i,j}+bx_i^2y_{i,j}+cx_iy_{i,j}^2+dx_i^2y_{i,j}^2\right)=0.
\end{split}
\end{equation}
At this point, one can see that $\BS h(\BS x)$ and $\BS H(\BS x)$ are
given by
\begin{equation}
\begin{split}
h_{i,j}(\BS x)&=\frac{\lambda_x}\varepsilon\left(ax_i+bx_i^2\right),\\
H_{i,j,i^\prime,j^\prime}(\BS y)&=\frac{\lambda_x}\varepsilon
\delta_{i,j}^{i^\prime,j^\prime}\left(cx_i+dx_i^2\right),
\end{split}
\end{equation}
in particular, $\BS H(\BS x)$ is a diagonal matrix. Also,
\begin{equation}
\frac 12\parderiv{^2 F_i}{\BS y^2}(\BS x):\BS\Sigma(\BS x)=
-\frac{\lambda_y}J\left(c+dx_i\right)\sum_{j=1}^J\Sigma_{i,j}^{i,j}(\BS x),
\end{equation}
that is, only the diagonal entries of the covariance matrix are
needed. This results in $\BS R^{h\to\bar z}$, $\BS R^{h\to\Sigma}$, $\BS
R^{H\to\bar z}$ and $\BS R^{H\to\Sigma}$ all being matrices rather than 3-
or 4-dimensional tensors.

\section{Numerical simulations}
\label{sec:num}

Here we show the results of numerical simulations with the new reduced
model for slow dynamics of the rescaled Lorenz model with nonlinear
coupling in \eqref{eq:lorenz_two_scale}. In particular compare the
numerical simulations for the three following systems:
\begin{enumerate}
\item The full two-scale rescaled Lorenz system from
  \eqref{eq:lorenz_two_scale};
\item The reduced model for slow dynamics from \eqref{eq:averaged_F}
  and \eqref{eq:z_Sigma};
\item A simplified version of \eqref{eq:averaged_F} with both the mean
  state $\bar{\BS z}$ and covariance matrix $\BS\Sigma$ of the fast
  variables are fixed at their average values $\bar{\BS z}^*$ and
  $\BS\Sigma^*$, without the correction terms from \eqref{eq:z_Sigma}
  (further referred to as the ``zero-order'' system).
\end{enumerate}
The quasi-Gaussian approximations in \eqref{eq:R_Gaussian} are used to
compute the approximations of the coupling terms in
\eqref{eq:z_Sigma}. For the time averaging in \eqref{eq:R_Gaussian} we
use long-term trajectories with averaging time window equals 10000
time units, while the correlation time window equals 50 time units (it
was observed that the time autocorrelation functions in
\eqref{eq:R_Gaussian} decays essentially to zero within the 50
time-unit window for all studied regimes). The reference mean state
$\BS{\bar z}^*$ and the covariance matrix $\BS\Sigma^*$ are also
computed by time-averaging of the full two-scale Lorenz model with the
same averaging window of 10000 time units.

The statistics of the model are invariant with respect to the index
permutation for the variables $x_i$, due to translational invariance
of the studied models. For the numerical study, we compute the
following long-term statistical quantities of $x_i$:
\begin{enumerate}[a.]
\item The probability density functions (PDF), computed by standard
  bin-counting.
\item The time autocorrelation functions $\langle
  x_i(t)x_i(t+s)\rangle$, where the angled brackets denote the time
  average is over $t$. These autocorrelation functions are normalized
  by the variance $\langle x_i^2\rangle$, so that the initial value at
  $s=0$ is always 1.
\item The time cross-correlation functions $\langle
  x_i(t)x_{i+1}(t+s)\rangle$, also normalized by the variance $\langle
  x_i^2\rangle$.
\item The energy autocorrelation function $$K(s)=\frac{\langle
  x_i^2(t)x_i^2(t+s)\rangle}{\langle x_i^2\rangle^2+2\langle
  x_i(t)x_i(t+s)\rangle^2}.$$ This energy autocorrelation function
  measures the non-Gaussianity of the process. It is identically 1 for
  all $s$ if the process is Gaussian (the Ornstein-Uhlenbeck process
  being an example). For details, see \cite{MajTimVan4}.
\end{enumerate}
The performance of the proposed approximation of the slow dynamics
depends on several factors. First, the precision will be affected by
the non-Gaussianity of the fast dynamics, since the quasi-Gaussian
linear response approximation is used for the computation of the
response operators. Second, the dependence of the mean state $\BS{\bar
  z}(\BS x)$ and covariance $\BS\Sigma(\BS x)$ for the fast variables
depends on the perturbations $\delta\BS h(\BS x)$ and $\delta\BS H(\BS
x)$ is generally nonlinear, and that should also affect the precision
of the approximation. Here we study the behavior of the proposed
approximation in variety of dynamical regimes of the rescaled Lorenz
model in \eqref{eq:lorenz_two_scale}. The following dynamical regimes
are studied:
\begin{itemize}
\item $N_x=20$, $J=4$ (so that $N_y=80$), so that the number of the
  fast variables is four times greater than the number of the slow
  variables.
\item $\varepsilon=0.01$. Typical geophysical processes, such as the
  annual and diurnal cycles, have the time scale separation of roughly
  two orders of magnitude.
\item $F_x=6$, $F_y=12$. The slow forcing $F_x$ adjusts the chaos and
  mixing properties of the slow variables, and in this work it is set
  to a weakly chaotic regime $F_x=6$. The fast forcing $F_y$ regulates
  chaos and mixing at the fast variables, which are usually more
  chaotic and mixing than the slow variables, so it is set to
  $F_y=12$.
\item $\lambda_x=\lambda_y=0.3$. This value of the coupling constant
  is chosen based on the previous work \cite{Abr9}, where the same
  value was used to test the method for linear coupling.
\item We test a mixture of different constants $(a,b,c,d)$ for each
  coupling term. First, we test the isolated coupling regimes with all
  constants but one set to zero, and then show the results for mixed
  coupled regimes, where different types of coupling are present
  simultaneously.
\end{itemize}

\subsection{Single type coupling}
\label{sec:single_coupling}

\begin{figure}
\picturehere{\figdir/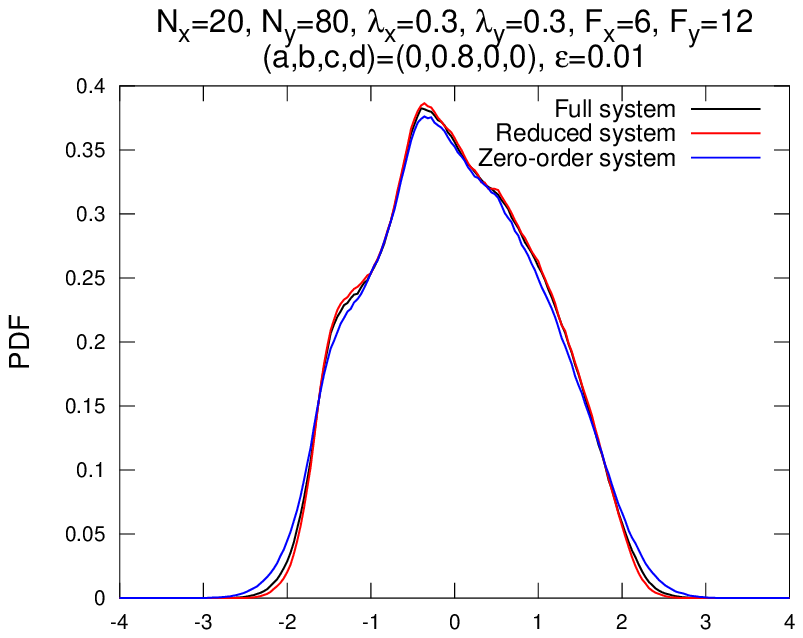}%
\picturehere{\figdir/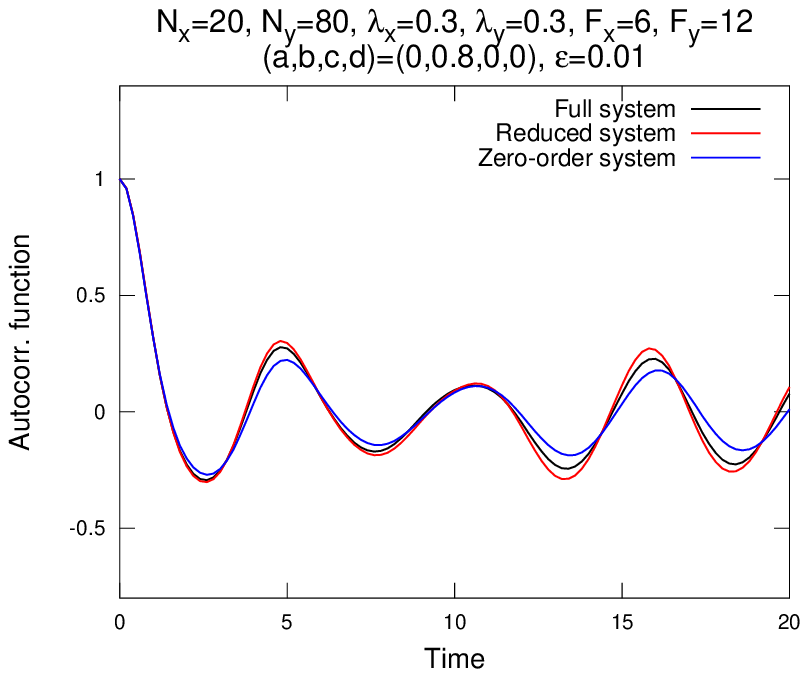}\\%
\picturehere{\figdir/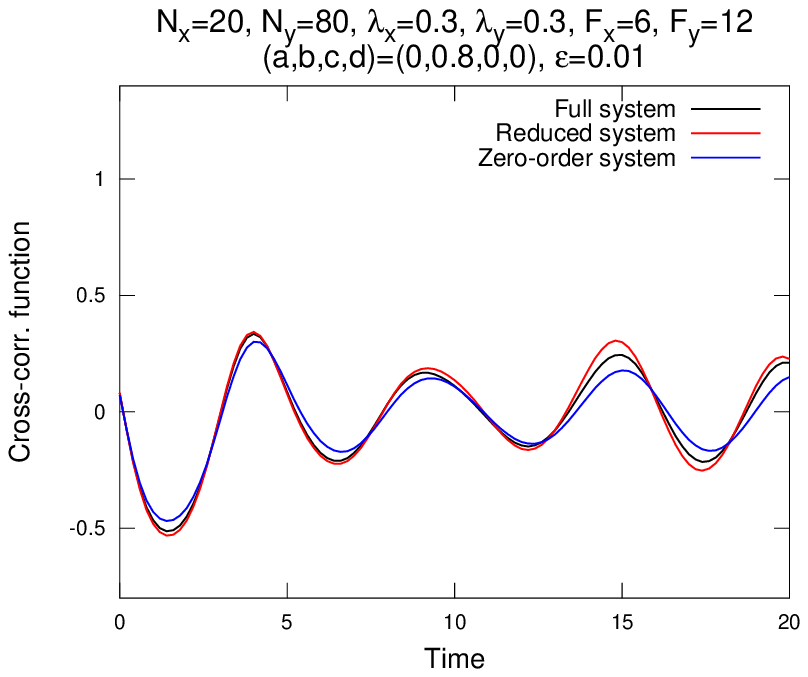}%
\picturehere{\figdir/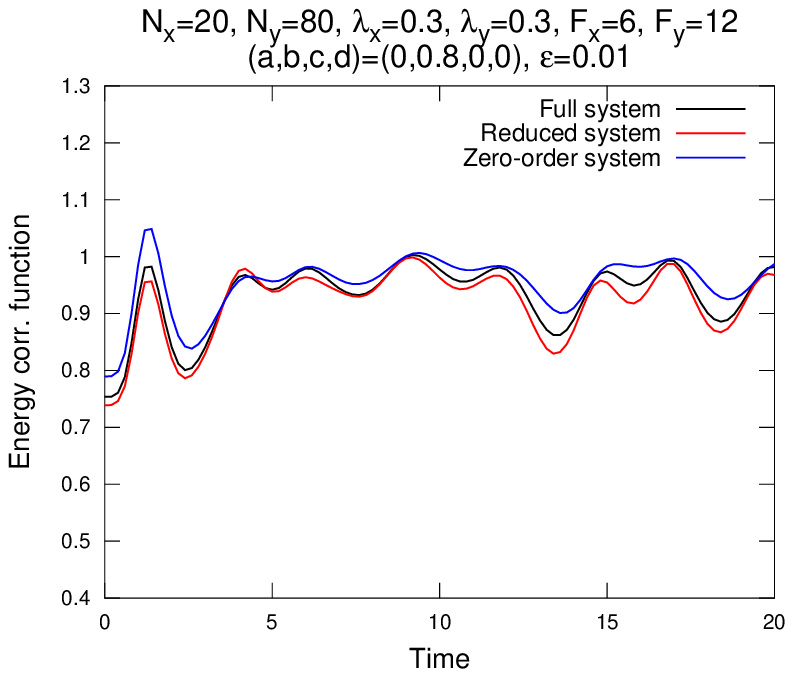}%
\caption{Probability density functions, time autocorrelations,
  cross-correlations and energy autocorrelations. Parameters:
  $(a,b,c,d)=(0,0.8,0,0)$.}
\label{fig:0_0.8_0_0}
\end{figure}
In this section we test the regimes with single type coupling (one of
the constants $a,b,c,d$ is nonzero, and the rest are zero). In this
section, we only test the regimes with nonzero constants $b,c,d$, as
the single type coupling regime with $a\neq 0$ corresponds to the
model previously studied in \cite{Abr9}.

\subsubsection{Regime with $(a,b,c,d)=(0,\pm 0.8,0,0)$}

\begin{figure}
\picturehere{\figdir/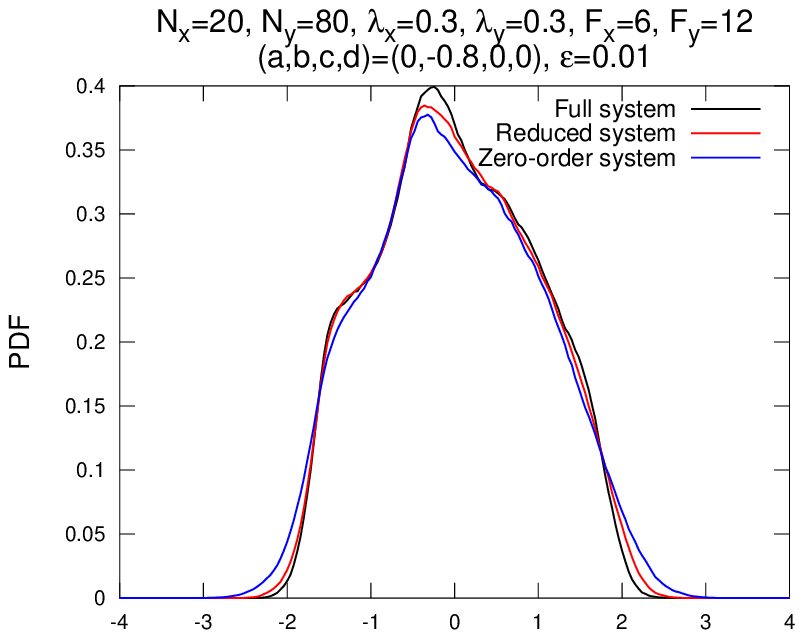}%
\picturehere{\figdir/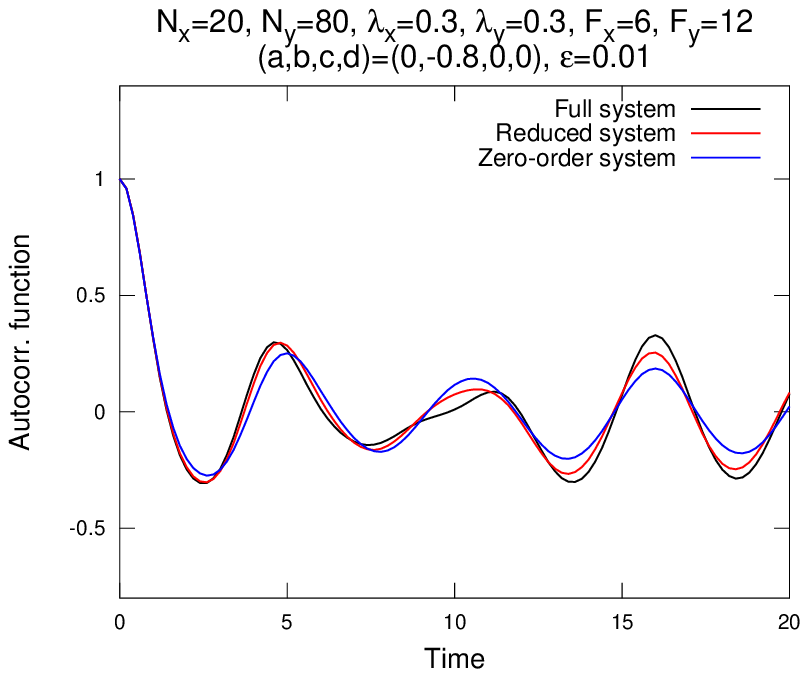}\\%
\picturehere{\figdir/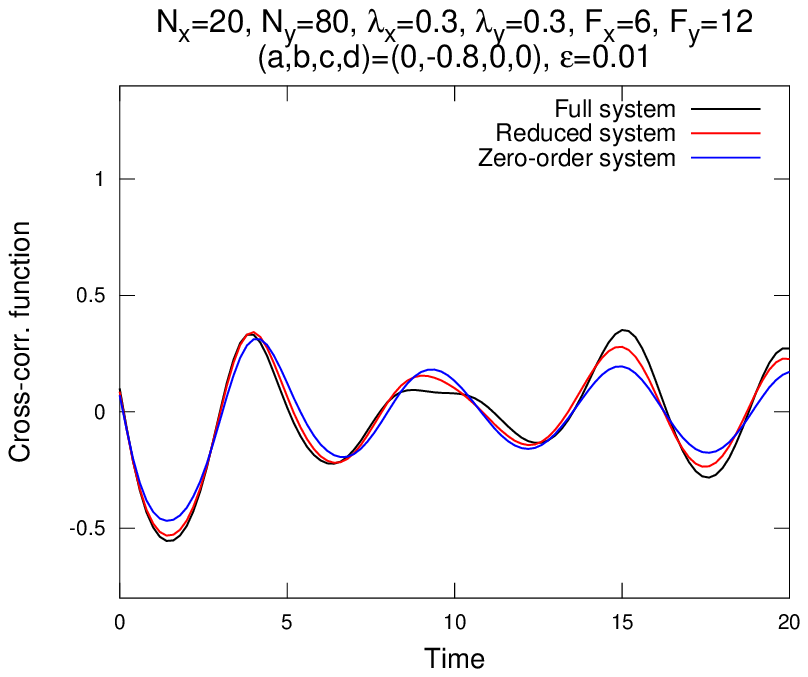}%
\picturehere{\figdir/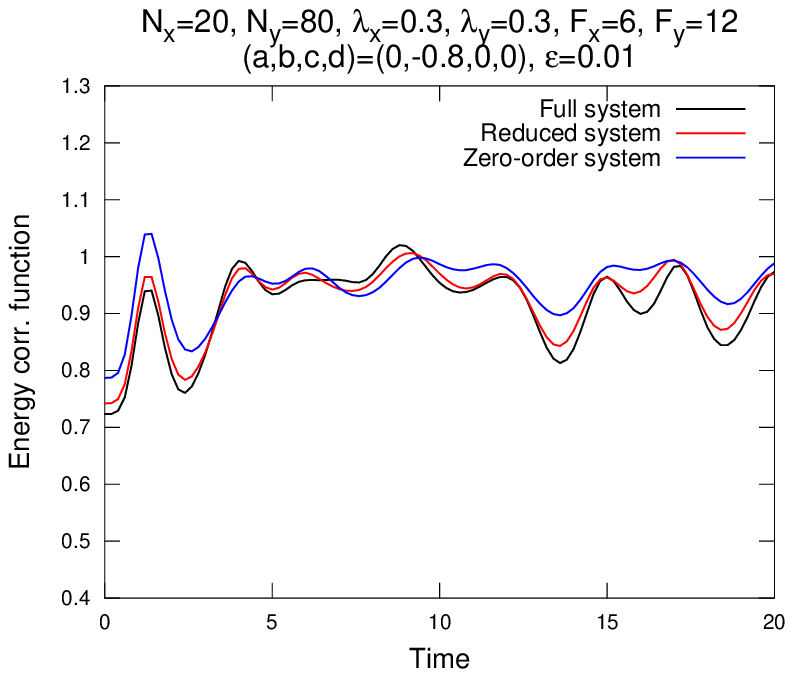}%
\caption{Probability density functions, time autocorrelations,
  cross-correlations and energy autocorrelations. Parameters:
  $(a,b,c,d)=(0,-0.8,0,0)$.}
\label{fig:0_-0.8_0_0}
\end{figure}
\begin{table}
\begin{tabular}{|c|}
\hline
$F_x=6$, $F_y=12$, $\lambda_x=\lambda_y=0.3$\\
\hline
\begin{tabular}{c|c}
$(a,b,c,d)=(0,0.8,0,0)$ & $(a,b,c,d)=(0,-0.8,0,0)$ \\
\hline\hline
\begin{tabular}{c|c|c}
 & Reduced & Zero-order \\
PDF & $3.103\cdot 10^{-3}$ & $6.61\cdot 10^{-3}$ \\
Corr. & $4.13\cdot 10^{-2}$ & $7.008\cdot 10^{-2}$ \\
C-corr. & $5.16\cdot 10^{-2}$ & $8.721\cdot 10^{-2}$ \\
K-corr. & $7.695\cdot 10^{-3}$ & $1.29\cdot 10^{-2}$ \\
\end{tabular}
&
\begin{tabular}{c|c|c}
 & Reduced & Zero-order \\
PDF & $6.471\cdot 10^{-3}$ & $1.473\cdot 10^{-2}$ \\
Corr. & $5.877\cdot 10^{-2}$ & $0.118$ \\
C-corr. & $7.009\cdot 10^{-2}$ & $0.1415$ \\
K-corr. & $9.033\cdot 10^{-3}$ & $2.384\cdot 10^{-2}$ \\
\end{tabular}
\end{tabular}\\
\hline
\end{tabular}
\caption{$L_2$-errors between the statistics of the slow variables of
  the full two-scale Lorenz model and the two reduced models, with
  $(a,b,c,d)=(0,\pm 0.8,0,0)$. Notations: ``Reduced'' stands for the
  reduced model from \eqref{eq:averaged_F} and \eqref{eq:z_Sigma}, and
  ``Zero-order'' stands for the poor man's version of the reduced
  model, with linear approximations for $\bar{\BS z}(\BS x)$ and
  $\BS\Sigma(\BS x)$ replaced by constant mean values $\bar{\BS z}^*$
  and $\BS\Sigma^*$.}
\label{tab:0_0.8_0_0}
\end{table}
Here we test the regimes with $b=\pm 0.8$, and $a=c=d=0$. This regime
corresponds to bilinear multiplicative coupling in the slow variables,
and quadratic additive coupling in the fast variables. The probability
density functions, time autocorrelation, cross-correlation, and energy
correlation functions are shown in Figures \ref{fig:0_0.8_0_0} and
\ref{fig:0_-0.8_0_0}, while the $L_2$-errors between the curves are
shown in Table \ref{tab:0_0.8_0_0}. There is little difference between
the plots, which means that the fast mean state $\BS{\bar z}(\BS x)$
and $\BS\Sigma(\BS x)$ are not very sensitive to changes in $\BS x$
for this type of coupling. Yet, one can see that the reduced model
with linear correction for $\BS{\bar z}(\BS x)$ and $\BS\Sigma(\BS x)$
more precisely captures statistics of the full scale model, than the
zero-order reduced model with fast mean state and covariance fixed at
$\BS{\bar z}^*$ and $\BS\Sigma^*$ (see Table \ref{tab:0_0.8_0_0}).
Additionally, it appears that there is a weak effect of chaos
suppression for this type of coupling regardless of the sign of $b$
(for both negative and positive $b$, the statistics of the coupled
model show somewhat slower decay of autocorrelations and
cross-correlations, and more sub-Gaussian energy correlation functions).

\subsubsection{Regime with $(a,b,c,d)=(0,0,\pm 0.3,0)$}

\begin{figure}
\picturehere{\figdir/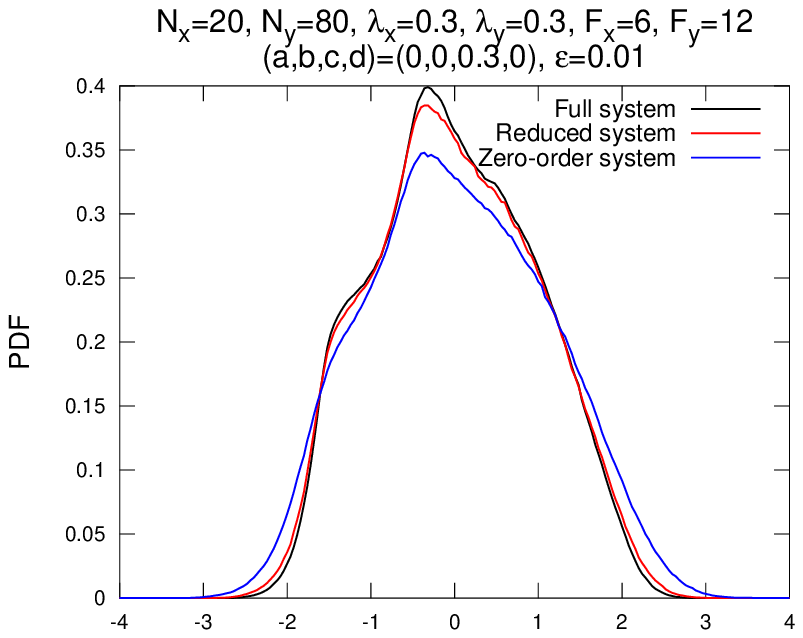}%
\picturehere{\figdir/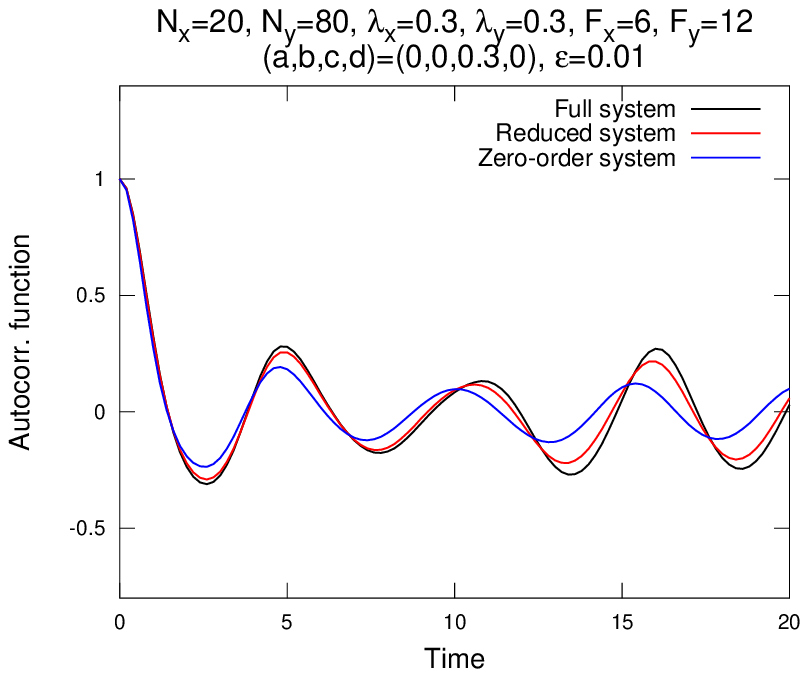}\\%
\picturehere{\figdir/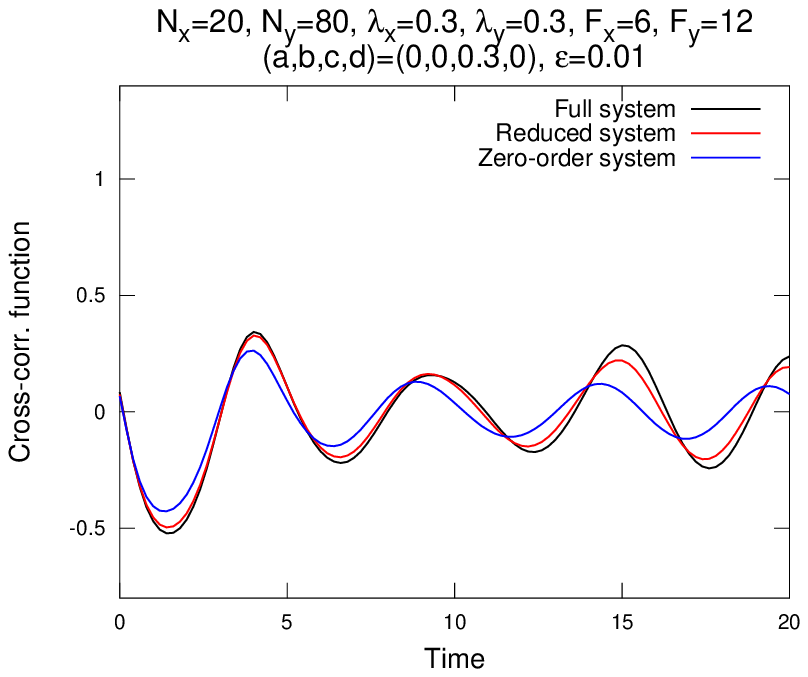}%
\picturehere{\figdir/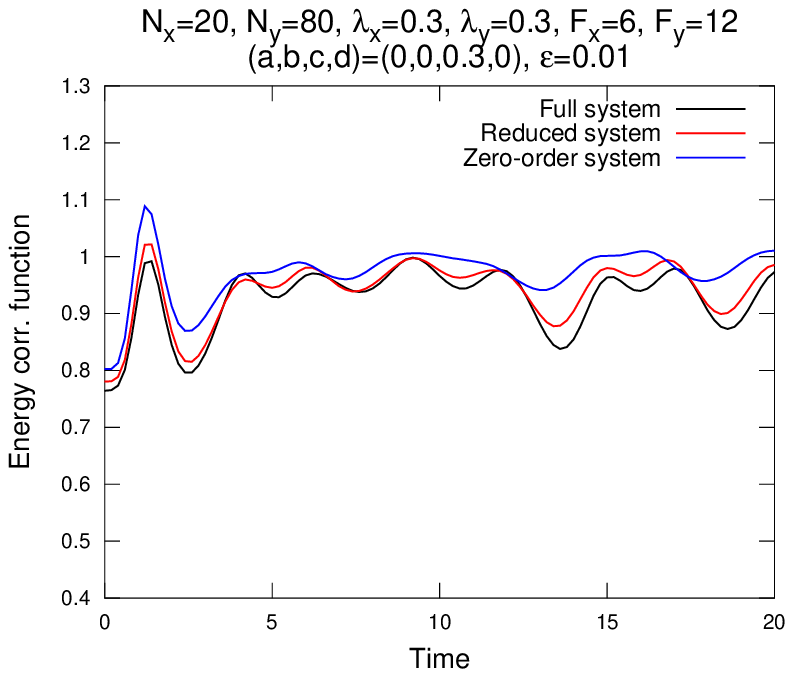}%
\caption{Probability density functions, time autocorrelations,
  cross-correlations and energy autocorrelations. Parameters:
  $(a,b,c,d)=(0,0,0.3,0)$.}
\label{fig:0_0_0.3_0}
\end{figure}
\begin{figure}
\picturehere{\figdir/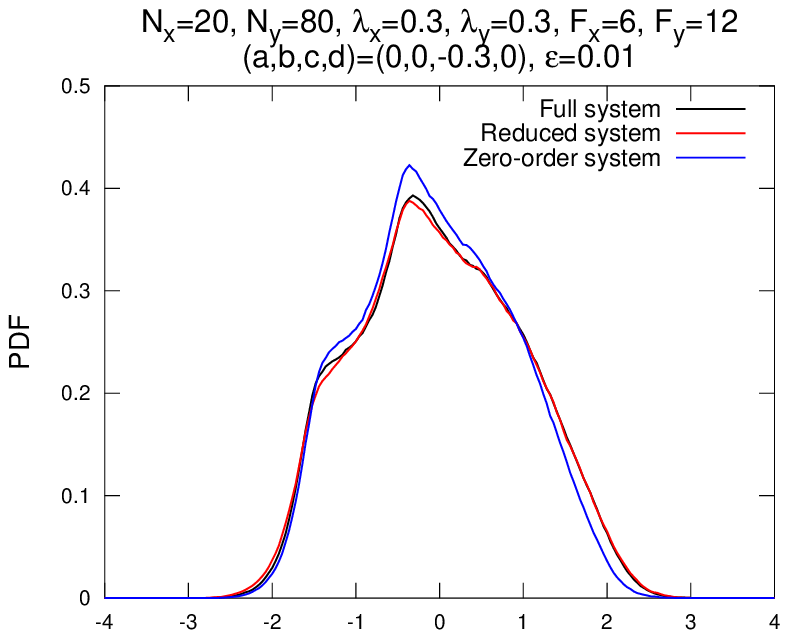}%
\picturehere{\figdir/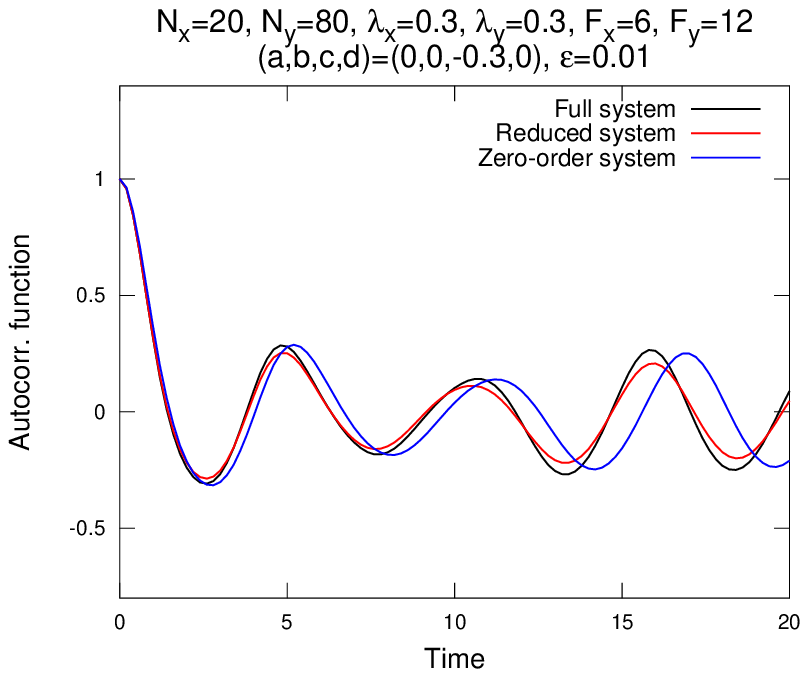}\\%
\picturehere{\figdir/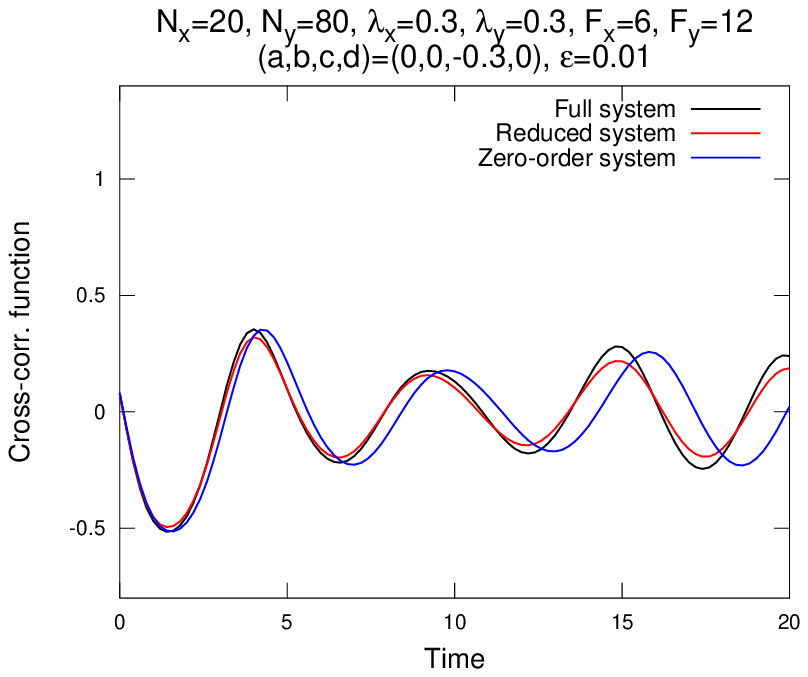}%
\picturehere{\figdir/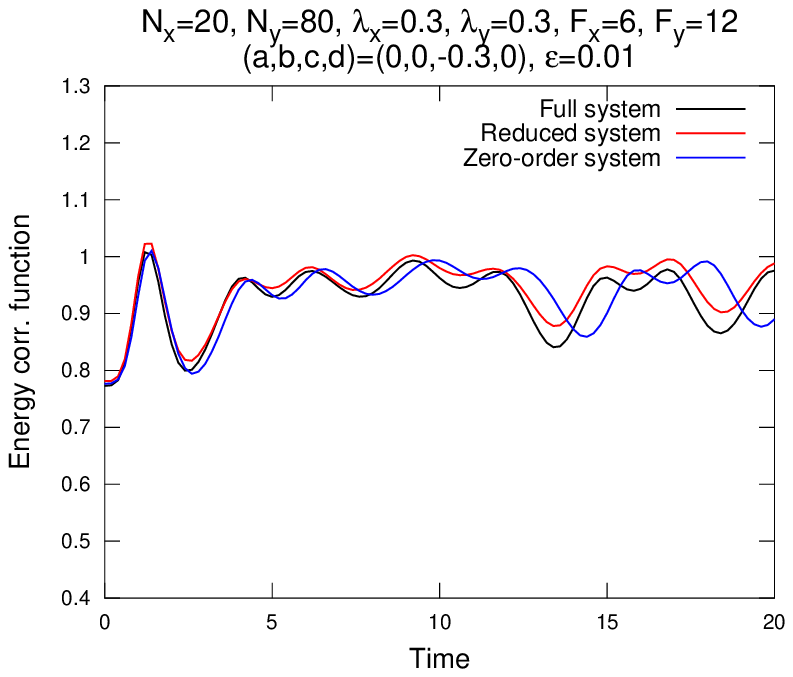}%
\caption{Probability density functions, time autocorrelations,
  cross-correlations and energy autocorrelations. Parameters:
  $(a,b,c,d)=(0,0,-0.3,0)$.}
\label{fig:0_0_-0.3_0}
\end{figure}
Here we test the regimes with $c=\pm 0.3$, and $a=b=d=0$. This regime
corresponds to quadratic additive coupling in the slow variables, and
bilinear multiplicative coupling in the fast variables. The
probability density functions, time autocorrelation,
cross-correlation, and energy correlation functions are shown in
Figures \ref{fig:0_0_0.3_0} and \ref{fig:0_0_-0.3_0}, while the
$L_2$-errors between the curves are shown in Table
\ref{tab:0_0_0.3_0}. Like in the previous case, there is little
difference between the plots, which means that the fast mean state
$\BS{\bar z}(\BS x)$ and $\BS\Sigma(\BS x)$ are not very sensitive to
changes in $\BS x$ for this type of coupling. Yet, one can see that
the reduced model with linear correction for $\BS{\bar z}(\BS x)$ and
$\BS\Sigma(\BS x)$ more precisely captures statistics of the full
scale model, than the zero-order reduced model with fast mean state
and covariance fixed at $\BS{\bar z}^*$ and $\BS\Sigma^*$ (see Table
\ref{tab:0_0_0.3_0}).  Additionally, it appears that there is a weak
effect of chaos suppression for this type of coupling for positive
value of $c=0.3$ (the statistics of the coupled model show somewhat
slower decay of autocorrelations and cross-correlations, and more
sub-Gaussian energy correlation). No such effect is observed for the
negative value $c=-0.3$.
\begin{table}
\begin{tabular}{|c|}
\hline
$F_x=6$, $F_y=12$, $\lambda_x=\lambda_y=0.3$\\
\hline
\begin{tabular}{c|c}
$(a,b,c,d)=(0,0,0.3,0)$ & $(a,b,c,d)=(0,0,-0.3,0)$ \\
\hline\hline
\begin{tabular}{c|c|c}
 & Reduced & Zero-order \\
PDF & $5.49\cdot 10^{-3}$ & $2.314\cdot 10^{-2}$ \\
Corr. & $5.559\cdot 10^{-2}$ & $0.1788$ \\
C-corr. & $6.265\cdot 10^{-2}$ & $0.2055$ \\
K-corr. & $9.605\cdot 10^{-3}$ & $2.659\cdot 10^{-2}$ \\
\end{tabular}
&
\begin{tabular}{c|c|c}
 & Reduced & Zero-order \\
PDF & $3.707\cdot 10^{-3}$ & $1.447\cdot 10^{-2}$ \\
Corr. & $5.238\cdot 10^{-2}$ & $0.2399$ \\
C-corr. & $6.271\cdot 10^{-2}$ & $0.2845$ \\
K-corr. & $9.62\cdot 10^{-3}$ & $2.041\cdot 10^{-2}$ \\
\end{tabular}
\end{tabular}\\
\hline
\end{tabular}
\caption{$L_2$-errors between the statistics of the slow variables of
  the full two-scale Lorenz model and the two reduced models, with
  $(a,b,c,d)=(0,0,\pm 0.3,0)$. Notations: ``Reduced'' stands for the
  reduced model from \eqref{eq:averaged_F} and \eqref{eq:z_Sigma}, and
  ``Zero-order'' stands for the poor man's version of the reduced
  model, with linear approximations for $\bar{\BS z}(\BS x)$ and
  $\BS\Sigma(\BS x)$ replaced by constant mean values $\bar{\BS z}^*$
  and $\BS\Sigma^*$.}
\label{tab:0_0_0.3_0}
\end{table}

\subsubsection{Regime with $(a,b,c,d)=(0,0,0,\pm 0.3)$}

\begin{figure}
\picturehere{\figdir/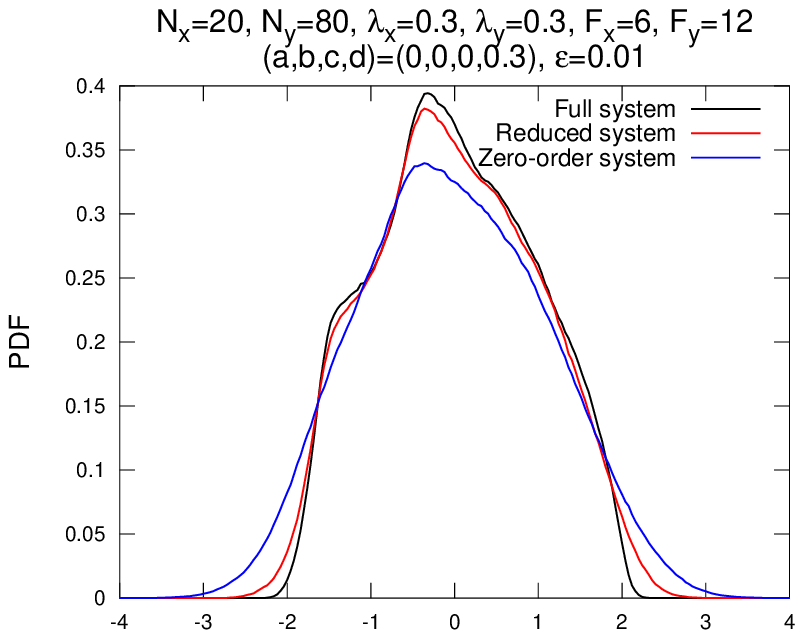}%
\picturehere{\figdir/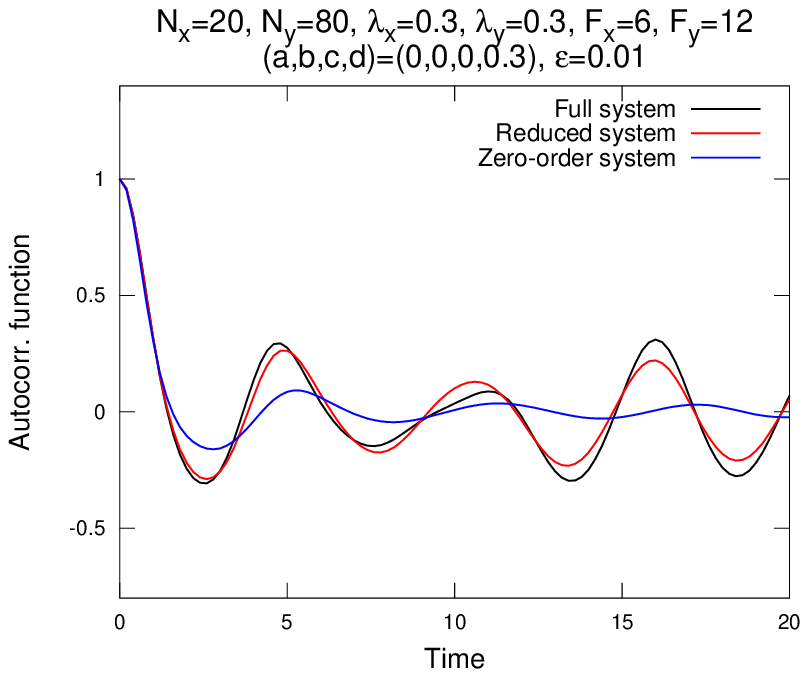}\\%
\picturehere{\figdir/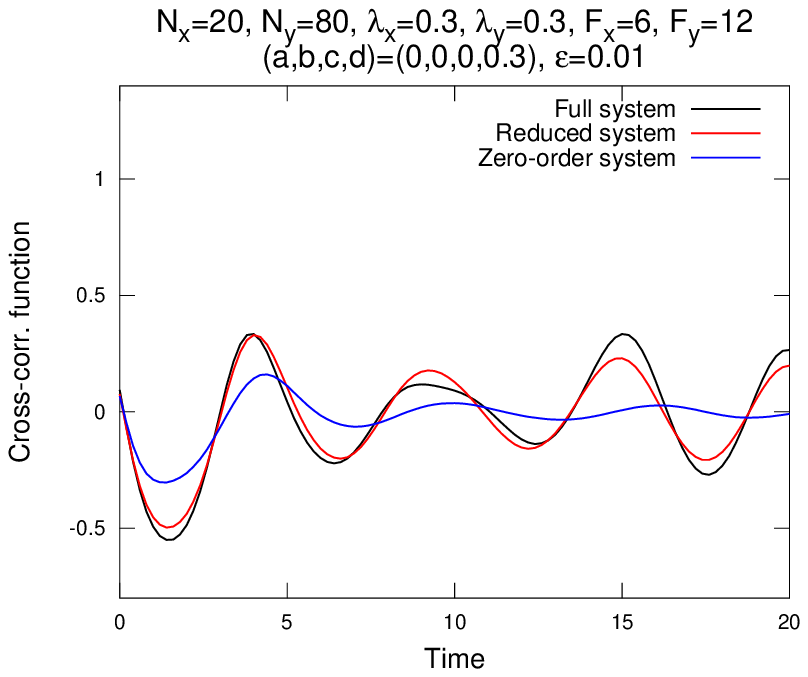}%
\picturehere{\figdir/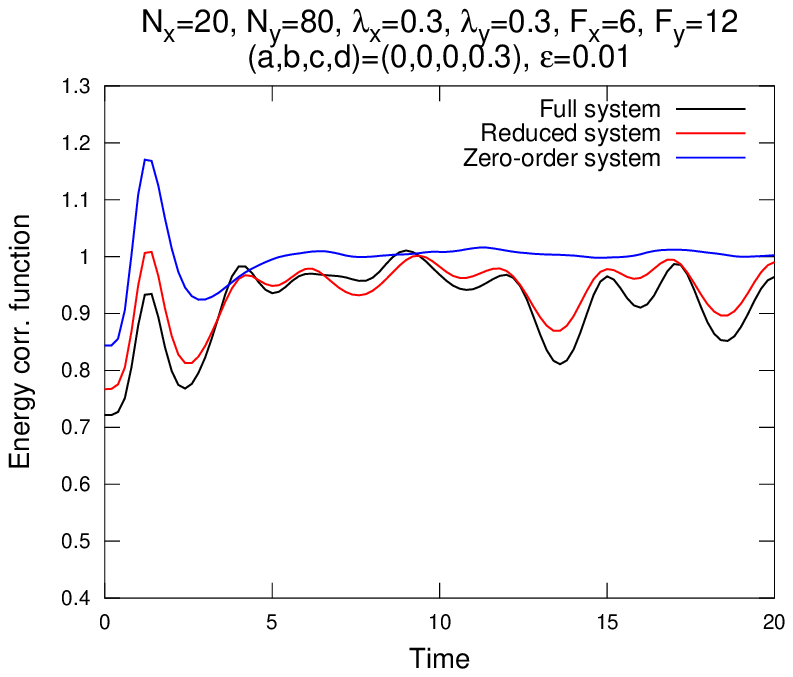}%
\caption{Probability density functions, time autocorrelations,
  cross-correlations and energy autocorrelations. Parameters:
  $(a,b,c,d)=(0,0,0,0.3)$.}
\label{fig:0_0_0_0.3}
\end{figure}
\begin{figure}
\picturehere{\figdir/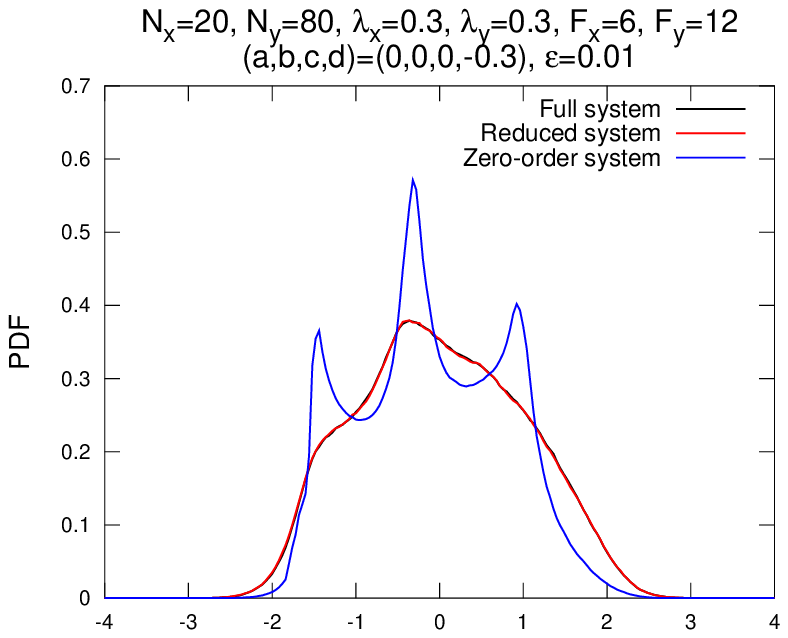}%
\picturehere{\figdir/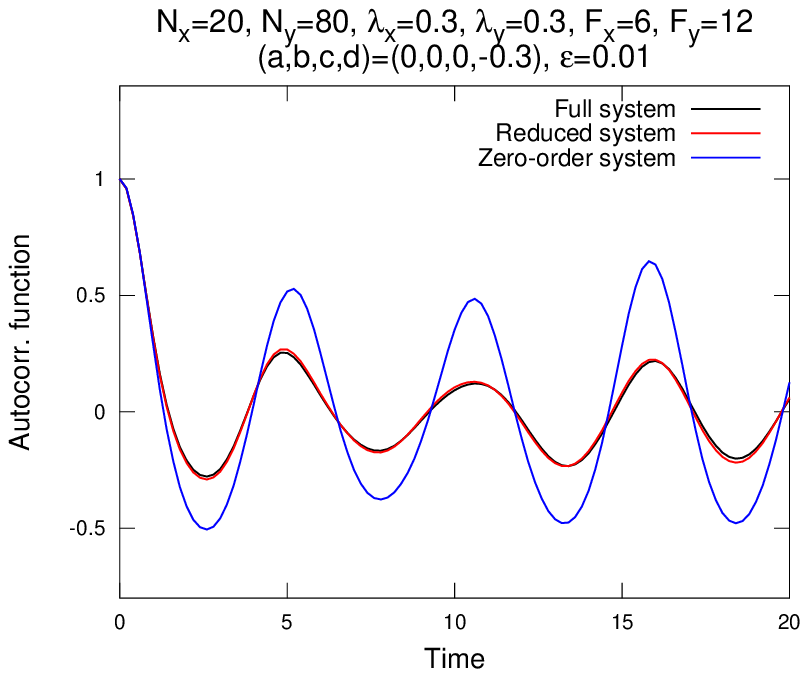}\\%
\picturehere{\figdir/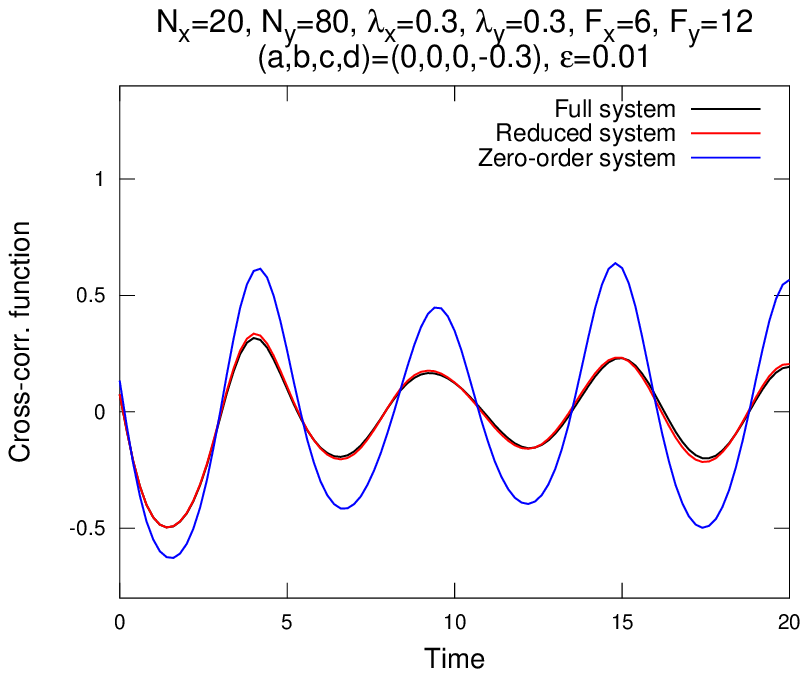}%
\picturehere{\figdir/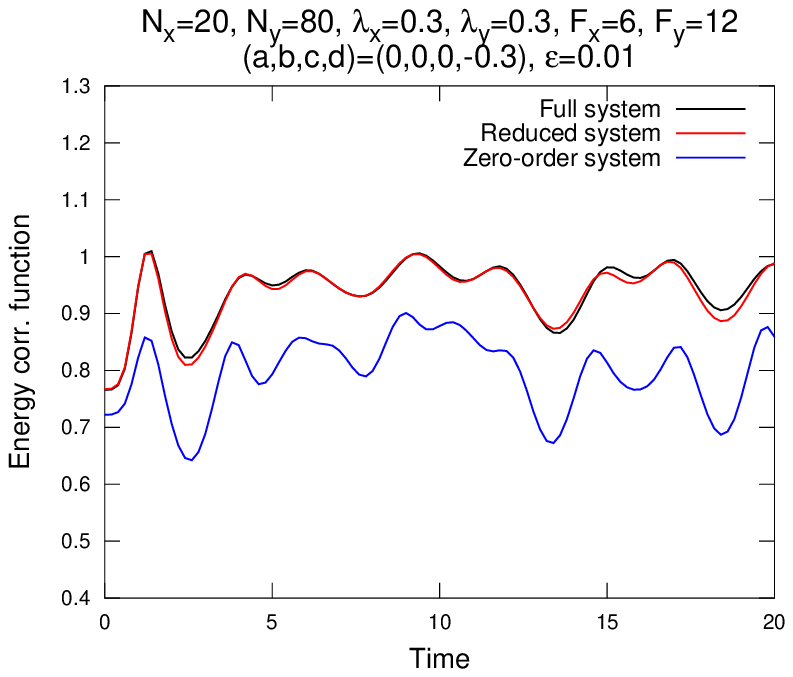}%
\caption{Probability density functions, time autocorrelations,
  cross-correlations and energy autocorrelations. Parameters:
  $(a,b,c,d)=(0,0,0,-0.3)$.}
\label{fig:0_0_0_-0.3}
\end{figure}
Here we test the regimes with $c=\pm 0.3$, and $a=b=d=0$. This regime
corresponds to quadratic multiplicative coupling in both the slow and
fast variables. The probability density functions, time
autocorrelation, cross-correlation, and energy correlation functions
are shown in Figures \ref{fig:0_0_0_0.3} and \ref{fig:0_0_0_-0.3},
while the $L_2$-errors between the curves are shown in Table
\ref{tab:0_0_0_0.3}. Unlike previous cases, there is a significant
difference between the results of the full two-scale/first-order
reduced model, and the zero-order reduced model, which means that the
fast mean state $\BS{\bar z}(\BS x)$ and $\BS\Sigma(\BS x)$ are
sensitive to changes in $\BS x$ for this type of coupling. Similar to
the previous cases, one can see that the reduced model with linear
correction for $\BS{\bar z}(\BS x)$ and $\BS\Sigma(\BS x)$ more
precisely captures the statistics of the full two-scale model, than
the zero-order reduced model with fast mean state and covariance fixed
at $\BS{\bar z}^*$ and $\BS\Sigma^*$ (see Table \ref{tab:0_0_0_0.3}).
Unlike previous types of coupling, here the chaos suppression or
amplification depends on the sign of the constant $d$. For the
positive value $d=0.3$, the coupled dynamics and the reduced model are
clearly less chaotic and mixing than the zero-order reduced model with
fast mean state and covariance fixed at $\BS{\bar z}^*$ and
$\BS\Sigma^*$, which follows from the difference in decay of the
correlation functions. The opposite effect is observed for the
negative value $d=-0.3$. In particular, observe that the PDF of the
zero-order reduced model displays three peaks for $d=-0.3$ (which is a
sign of quasi-periodic motion), while the PDFs of both the two-scale
Lorenz model and reduced models are unimodal.
\begin{table}
\begin{tabular}{|c|}
\hline
$F_x=6$, $F_y=12$, $\lambda_x=\lambda_y=0.3$\\
\hline
\begin{tabular}{c|c}
$(a,b,c,d)=(0,0,0,0.3)$ & $(a,b,c,d)=(0,0,0,-0.3)$ \\
\hline\hline
\begin{tabular}{c|c|c}
 & Reduced & Zero-order \\
PDF & $1.012\cdot 10^{-2}$ & $2.935\cdot 10^{-2}$ \\
Corr. & $7.397\cdot 10^{-2}$ & $0.2638$ \\
C-corr. & $9.183\cdot 10^{-2}$ & $0.3196$ \\
K-corr. & $1.632\cdot 10^{-2}$ & $5.002\cdot 10^{-2}$ \\
\end{tabular}
&
\begin{tabular}{c|c|c}
 & Reduced & Zero-order \\
PDF & $1.342\cdot 10^{-3}$ & $5.272\cdot 10^{-2}$ \\
Corr. & $1.923\cdot 10^{-2}$ & $0.3881$ \\
C-corr. & $2.497\cdot 10^{-2}$ & $0.4624$ \\
K-corr. & $3.656\cdot 10^{-3}$ & $6.977\cdot 10^{-2}$ \\
\end{tabular}
\end{tabular}\\
\hline
\end{tabular}
\caption{$L_2$-errors between the statistics of the slow variables of
  the full two-scale Lorenz model and the two reduced models, with
  $(a,b,c,d)=(0,0,0,\pm 0.3)$. Notations: ``Reduced'' stands for the
  reduced model from \eqref{eq:averaged_F} and \eqref{eq:z_Sigma}, and
  ``Zero-order'' stands for the poor man's version of the reduced
  model, with linear approximations for $\bar{\BS z}(\BS x)$ and
  $\BS\Sigma(\BS x)$ replaced by constant mean values $\bar{\BS z}^*$
  and $\BS\Sigma^*$.}
\label{tab:0_0_0_0.3}
\end{table}
\begin{figure}
\picturehere{\figdir/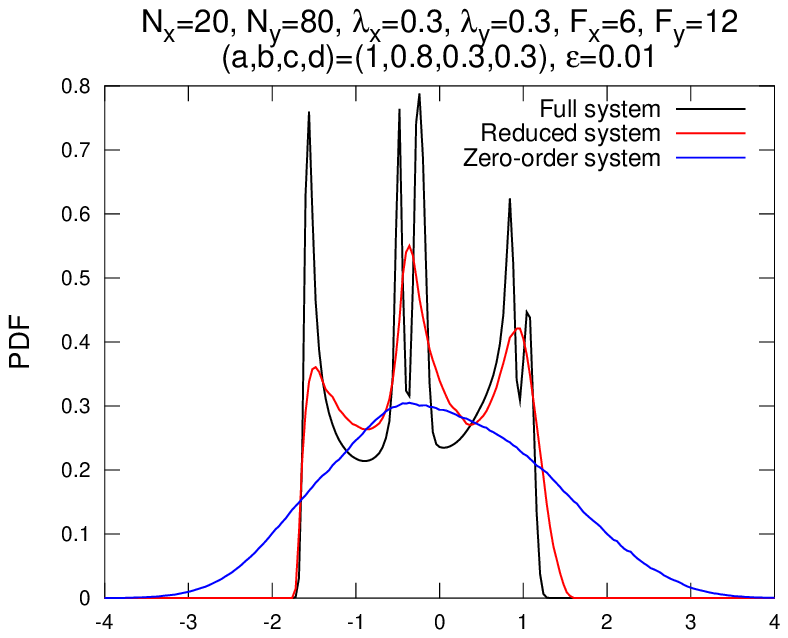}%
\picturehere{\figdir/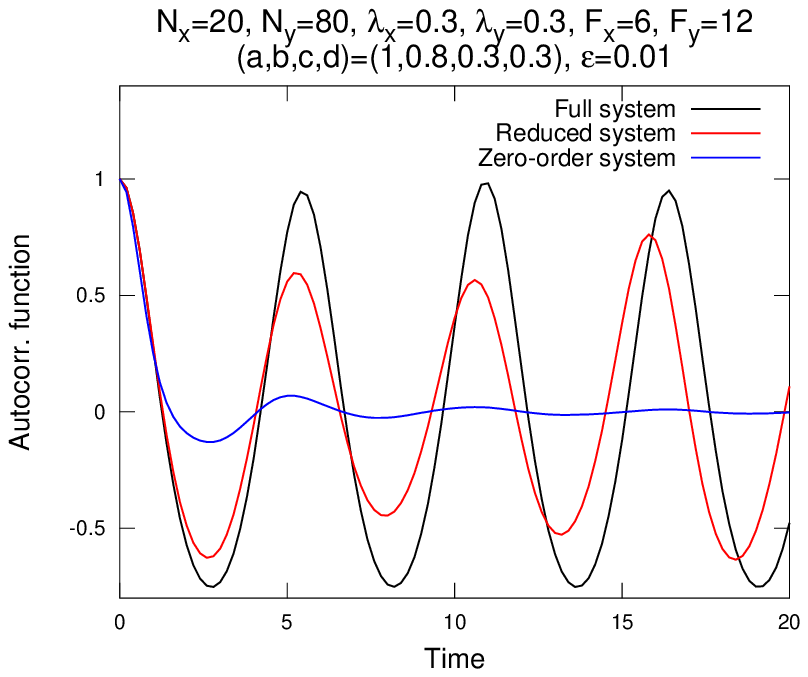}\\%
\picturehere{\figdir/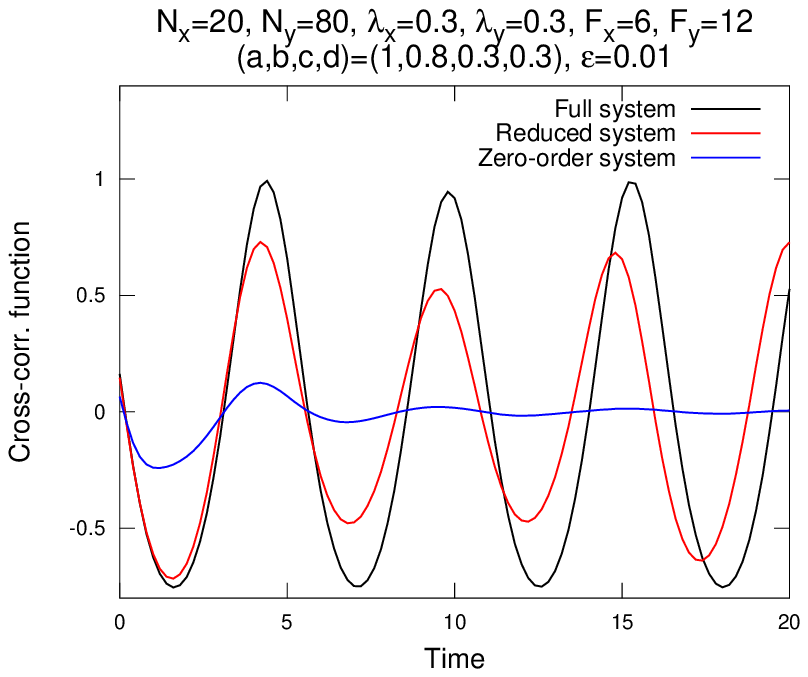}%
\picturehere{\figdir/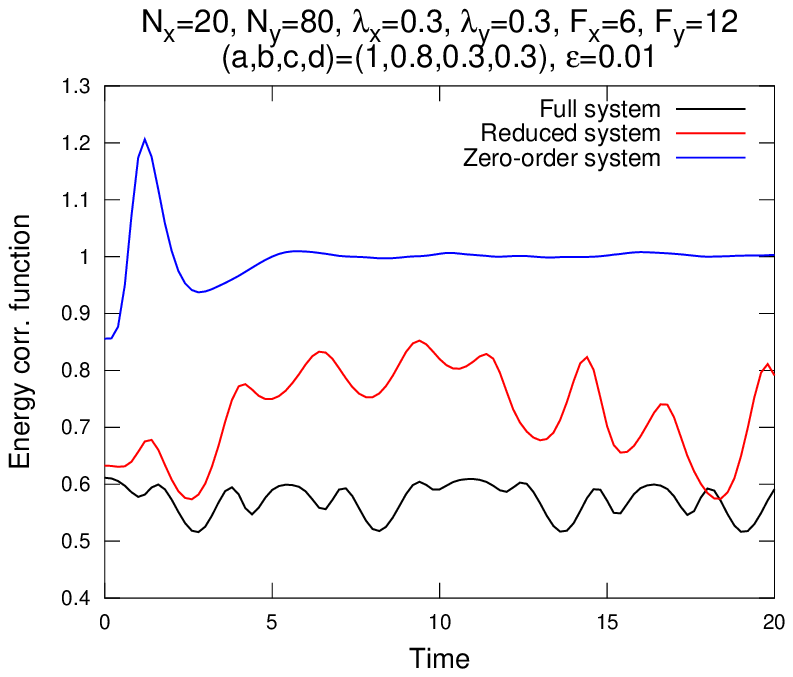}%
\caption{Probability density functions, time autocorrelations,
  cross-correlations and energy autocorrelations. Parameters:
  $(a,b,c,d)=(1,0.8,0.3,0.3)$.}
\label{fig:1_0.8_0.3_0.3}
\end{figure}

\subsection{Combined coupling}
\label{sec:combined_coupling}

\begin{figure}
\picturehere{\figdir/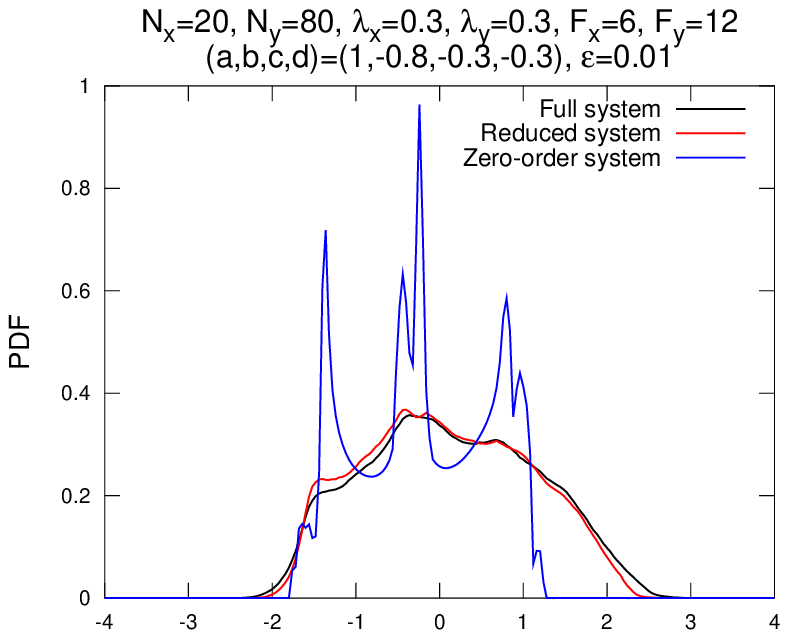}%
\picturehere{\figdir/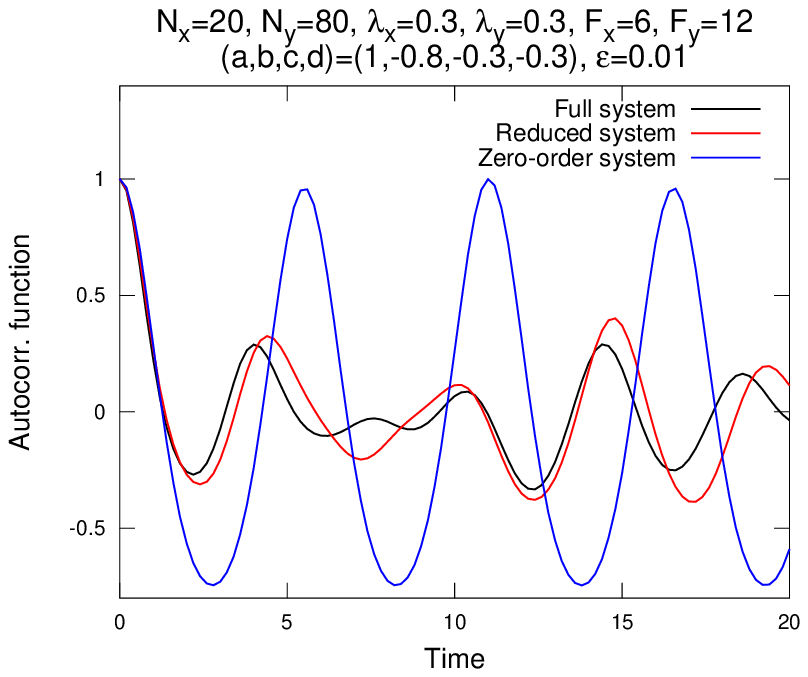}\\%
\picturehere{\figdir/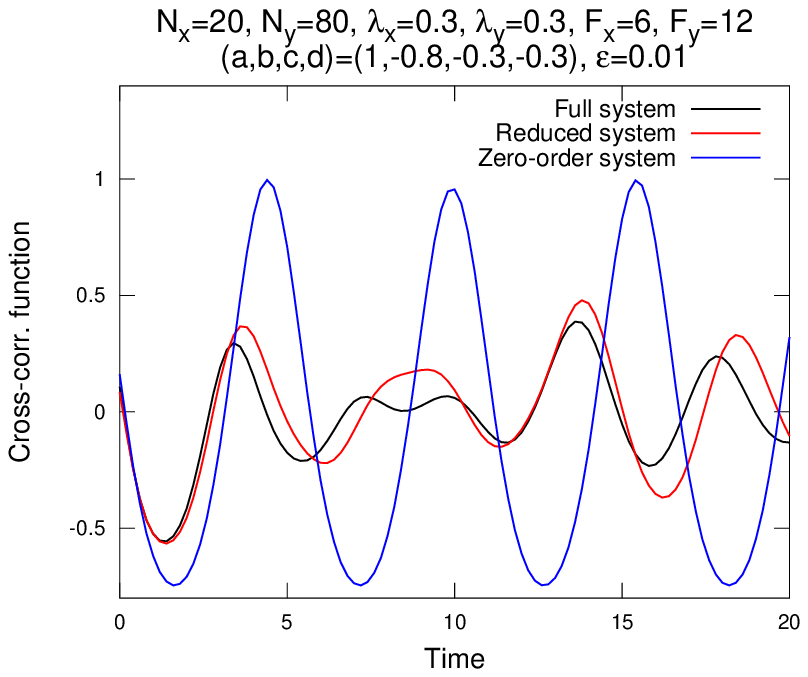}%
\picturehere{\figdir/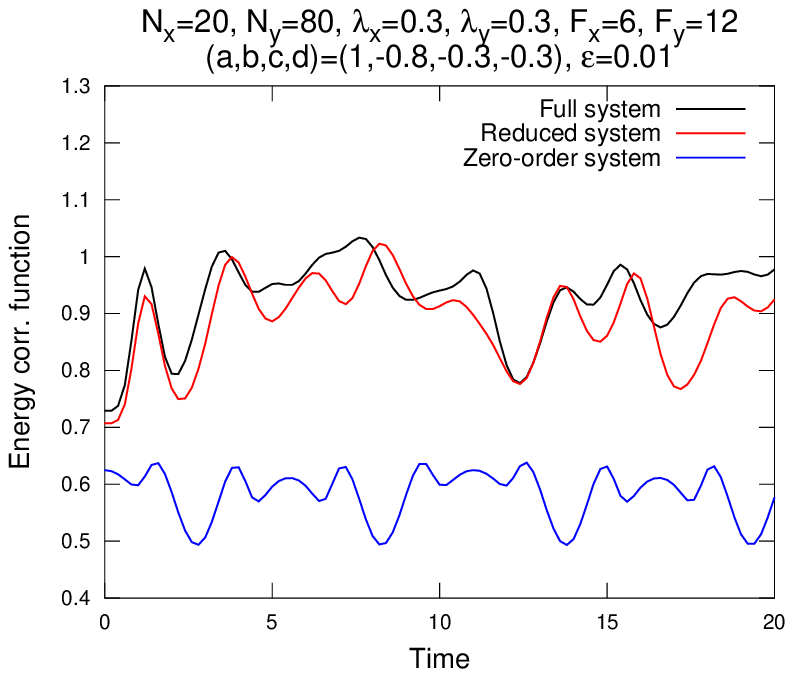}%
\caption{Probability density functions, time autocorrelations,
  cross-correlations and energy autocorrelations. Parameters:
  $(a,b,c,d)=(1,-0.8,-0.3,-0.3)$.}
\label{fig:1_-0.8_-0.3_-0.3}
\end{figure}
In this section we test the regimes with all types of coupling
observed previously in Section \ref{sec:single_coupling}, combined
together in the Lorenz model. In addition, the linear coupling is also
switched on by setting $a=1$. Here we study the sets of parameters
$(a,b,c,d)=(1,\pm 0.8,\pm 0.3,\pm 0.3)$. The probability density
functions, time autocorrelation, cross-correlation, and energy
correlation functions are shown in Figures \ref{fig:1_0.8_0.3_0.3} and
\ref{fig:1_-0.8_-0.3_-0.3}, while the $L_2$-errors between the curves
are shown in Table \ref{tab:1_0.8_0.3_0.3}. It turns out that when all
types of coupling are combined together there is a major difference
between the two-scale/first-order reduced models, and the zero-order
reduced model with fast mean state and covariance fixed at $\BS{\bar
  z}^*$ and $\BS\Sigma^*$, which means that the fast mean state
$\BS{\bar z}(\BS x)$ and $\BS\Sigma(\BS x)$ are quite sensitive to
changes in $\BS x$. Similar to the previous cases, here one can see
that the reduced model with linear correction for $\BS{\bar z}(\BS x)$
and $\BS\Sigma(\BS x)$ much more precisely captures statistics of the
full two-scale scale model, than the zero-order reduced model (see
Table \ref{tab:0_0_0_0.3}). Here the chaos suppression or
amplification depends on the signs of the constant $b$, $c$ and $d$.
For positive values of these constants, the coupled dynamics and the
reduced model are clearly less chaotic and mixing than the zero-order
reduced model, which follows from the difference in decay of the
correlation functions. The opposite effect is observed for negative
values of $b$, $c$ and $d$. In particular, observe that the PDF of the
zero-order reduced model displays three strong peaks for $d=-0.3$
(which is a sign of quasi-periodic motion), while the PDFs of both the
two-scale Lorenz model and reduced models are unimodal. On the
contrary, for $d=0.3$ the PDF of the two-scale Lorenz model displays
three peaks, which are roughly captured by the reduced model, while the
zero-order reduced model produces nearly Gaussian PDF.
\begin{table}
\begin{tabular}{|c|}
\hline
$F_x=6$, $F_y=12$, $\lambda_x=\lambda_y=0.3$\\
\hline
\begin{tabular}{c|c}
$(a,b,c,d)=(1,0.8,0.3,0.3)$ & $(a,b,c,d)=(1,-0.8,-0.3,-0.3)$ \\
\hline\hline
\begin{tabular}{c|c|c}
 & Reduced & Zero-order \\
PDF & $6.837\cdot 10^{-2}$ & $0.1121$ \\
Corr. & $0.2237$ & $0.4083$ \\
C-corr. & $0.2322$ & $0.4153$ \\
K-corr. & $0.134$ & $0.3356$ \\
\end{tabular}
&
\begin{tabular}{c|c|c}
 & Reduced & Zero-order \\
PDF & $1.24\cdot 10^{-2}$ & $0.1252$ \\
Corr. & $0.2181$ & $1.209$ \\
C-corr. & $0.2477$ & $1.337$ \\
K-corr. & $2.89\cdot 10^{-2}$ & $0.1707$ \\
\end{tabular}
\end{tabular}\\
\hline
\end{tabular}
\caption{$L_2$-errors between the statistics of the slow variables of
  the full two-scale Lorenz model and the two reduced models, with
  $(a,b,c,d)=(1,\pm 0.8,\pm 0.3,\pm 0.3)$. Notations: ``Reduced''
  stands for the reduced model from \eqref{eq:averaged_F} and
  \eqref{eq:z_Sigma}, and ``Zero-order'' stands for the poor man's
  version of the reduced model, with linear approximations for
  $\bar{\BS z}(\BS x)$ and $\BS\Sigma(\BS x)$ replaced by constant
  mean values $\bar{\BS z}^*$ and $\BS\Sigma^*$.}
\label{tab:1_0.8_0.3_0.3}
\end{table}

\section{Conclusions}
\label{sec:sum}

In this work we develop a simple approach for approximation of
multiscale dynamics with nonlinear and multiplicative coupling via a
reduced model for slow variables alone. The method is based on the
linear approximation of averaged coupling terms by means of the
fluctuation-dissipation theorem, which is only requires a single
computation of certain long-term statistics from the fast limiting
system with the slow terms set as constant parameters. This work is a
direct extension of \cite{Abr9} onto nonlinear and multiplicative
coupling (the original method in \cite{Abr9} was developed for linear
coupling in both slow and fast variables). We verify through the
numerical simulations with the rescaled two-scale Lorenz 96 model
\cite{Abr8,Abr9} that, with nonlinear and multiplicative coupling in
both slow and fast variables, the new simple reduced model produces
statistics which are consistent with those of the complete two-scale
Lorenz model. In contrast, the ``zero-order'' reduced model with
constant parameterization of fast variables in coupling terms fails to
reproduce the same set of statistics with comparable precision. The
method appears to be convenient for practical applications due to its
explicit construction -- it lacks unknown parameters which have to be
determined implicitly by comparing the performance of the reduced
model against the full multiscale dynamics.

\appendix
\section{Detailed derivation of the closure terms for the mean state
and covariance matrix}
\label{sec:app_reduced}

In order to derive the formulas for the linear responses of the mean
and covariance of the fast variables, for convenience we first make a
linear change of variables in the following way: we rewrite
\eqref{eq:dyn_sys_fast_limiting_delta} as
\begin{equation}
\begin{split}
\deriv{\BS q}\tau=\BS S^{-1}\BS g(\BS{\bar z}^*+\BS S\BS q)+\BS
S^{-1}(\BS H^*+\delta\BS H(\BS x))\BS S\BS q+\\+\BS S^{-1}(\BS
h^*+\BS H^*\BS{\bar z}^*+\delta\BS h(\BS x)+\delta\BS H(\BS
x)\BS{\bar z}^*),
\end{split}
\end{equation}
where $\BS S=\BS\Sigma^{*\frac 12}$, and $\BS z=\BS{\bar z}^*+\BS S\BS q$,
that is, $\BS q$ is the fluctuation of $\BS z$ around $\BS{\bar z}^*$
with the identity covariance matrix. Changing notations, we obtain
\begin{equation}
\label{eq:dyn_sys_fast_limiting_*}
\begin{split}
\deriv{\BS q}{\BS\tau}=\BS S^{-1}\BS g&(\BS{\bar z}^*+\BS S\BS q)
+(\BS{\hat H}^*+\delta\BS{\hat H}(\BS x))\BS q+\BS{\hat h}^*+
\delta\BS{\hat h}(\BS x),\\
\BS{\hat h}^*&=\BS S^{-1}(\BS h^*+\BS H^*\BS{\bar z}^*),\quad
\BS{\hat H}^*=\BS S^{-1}\BS H^*\BS S,\\
\delta\BS{\hat h}(\BS x)=\BS S^{-1}&(\delta\BS h(\BS x)+
\delta\BS H(\BS x)\BS{\bar z}^*),\quad
\delta\BS{\hat H}(\BS x)=\BS S^{-1}\delta\BS H(\BS x)\BS S.
\end{split}
\end{equation}
Here, we can consider \eqref{eq:dyn_sys_fast_limiting_*} as a
dynamical system perturbed by $\delta\BS{\hat h}$ and $\delta\BS{\hat
  H}$, which has zero mean state and identity covariance matrix in the
unperturbed state, and use the FDT to estimate the linear response of
the mean state $\langle\BS q\rangle$ and the covariance $\langle\BS
q\BS q^T\rangle$ to these perturbations, which can later be mapped
into the original coordinates via backward linear transformation. The
linear responses of the mean state $\delta\langle\BS q\rangle$ and
covariance $\delta\langle\BS q\BS q^T\rangle$ are given by,
respectively,
\begin{equation}
\begin{split}
\delta\langle\BS q\rangle_i=\BS R_{ij}^{\hat h\to\hat{\bar z}}\delta\BS{\hat
  h}_j +\BS R_{ijk}^{\hat H\to\hat{\bar z}}\delta\BS{\hat
  H}_{jk},\\ \delta\langle\BS q\BS q^T\rangle_{ij}=\BS R_{ijk}^{\hat
  h\to\hat\Sigma} \delta\BS{\hat h}_k+\BS R_{ijkl}^{\hat
  H\to\hat\Sigma}\delta\BS{\hat H}_{kl},
\end{split}
\end{equation}
where the linear response operators are given by
\begin{subequations}
\begin{equation}
\BS R_{ij}^{\hat h\to\hat{\bar z}}=\int_0^\infty\int_{\mathbb R^{N_y}}
\parderiv{}{q_j}(\phi^s\BS q)_i\dif\mu(\BS q)\dif s,
\end{equation}
\begin{equation}
\BS R_{ijk}^{\hat H\to\hat{\bar z}}=\int_0^\infty\int_{\mathbb R^{N_y}}
\parderiv{}{q_j}(\phi^s\BS q)_iq_k\dif\mu(\BS q)\dif s,
\end{equation}
\begin{equation}
\BS R_{ijk}^{\hat h\to\hat\Sigma}=\int_0^\infty\int_{\mathbb R^{N_y}}
\parderiv{}{q_k}(\phi^s\BS q\otimes\phi^s\BS q)_{ij}
\dif\mu(\BS q)\dif s,
\end{equation}
\begin{equation}
\BS R_{ijkl}^{\hat H\to\hat\Sigma}=\int_0^\infty\int_{\mathbb R^{N_y}}
\parderiv{}{q_k}(\phi^s\BS q\otimes\phi^s\BS q)_{ij}q_l
\dif\mu(\BS q)\dif s,
\end{equation}
\end{subequations}
where $\phi^s$ is the flow generated by the unperturbed system in
\eqref{eq:dyn_sys_fast_limiting_*} with $\delta\hat{\BS h}$ and
$\delta\hat{\BS H}$ set to zeros, and $\dif\mu$ is its invariant
measure \cite{Abr5,Abr6,Abr7,Abr8,Abr9,AbrMaj4,AbrMaj5,AbrMaj6}. At
this point, we are going to assume that $\mu$ has a Gaussian density
with zero mean state and identity covariance matrix (so-called
quasi-Gaussian FDT approximation,
\cite{Abr9,AbrMaj4,AbrMaj5,AbrMaj6,MajAbrGro}). Then, after
integrating by parts, the responses of the mean state and covariance
can be computed as
\begin{subequations}
\begin{equation}
\BS R_{ij}^{\hat h\to\hat{\bar z}}=\int_0^\infty\int_{\mathbb R^{N_y}}
(\phi^s\BS q)_iq_j\dif\mu(\BS q)\dif s,
\end{equation}
\begin{equation}
\BS R_{ijk}^{\hat H\to\hat{\bar z}}=\int_0^\infty\int_{\mathbb R^{N_y}}
(\phi^s\BS q)_iq_jq_k\dif\mu(\BS z)\dif s,
\end{equation}
\begin{equation}
\BS R_{ijk}^{\hat h\to\hat\Sigma}=\int_0^\infty\int_{\mathbb R^{N_y}}
(\phi^s\BS q\otimes\phi^s\BS q)_{ij}q_k\dif\mu(\BS z)
\dif s,
\end{equation}
\begin{equation}
\BS R_{ijkl}^{\hat H\to\hat\Sigma}=\int_0^\infty\left[\int_{\mathbb R^{N_y}}
(\phi^s\BS q\otimes\phi^s\BS q)_{ij}q_kq_l
\dif\mu(\BS z)-\delta_{ij}\delta_{kl}\right]\dif s.
\end{equation}
\end{subequations}
By replacing measure averages with time averages (under the ergodicity
assumption), we obtain
\begin{subequations}
\begin{equation}
\BS R_{ij}^{\hat h\to\hat{\bar z}}=\int_0^\infty\left[\lim_{r\to\infty}
\frac 1r\int_0^rq_i(t+s)q_j(t)\dif t\right]\dif s,
\end{equation}
\begin{equation}
\BS R_{ijk}^{\hat H\to\hat{\bar z}}=\int_0^\infty\left[\lim_{r\to\infty}\frac 1r
\int_0^rq_i(t+s)q_j(t)q_k(t)\dif t\right]\dif s,
\end{equation}
\begin{equation}
\BS R_{ijk}^{\hat h\to\hat\Sigma}=\int_0^\infty\left[\lim_{r\to\infty}
\frac 1r\int_0^rq_i(t+s)q_j(t+s)q_k(t)\dif t\right]\dif s,
\end{equation}
\begin{equation}
\BS R_{ijkl}^{\hat H\to\hat\Sigma}=\int_0^\infty\left[\lim_{r\to\infty}
\frac 1r\int_0^rq_i(t+s)q_j(t+s)q_k(t)q_l(t)\dif t-\delta_{ij}
\delta_{kl}\right]\dif s.
\end{equation}
\end{subequations}
Returning back to the original response coordinates, we obtain
\begin{subequations}
\begin{equation}
\BS R_{ij}^{\hat h\to\bar z}=\int_0^\infty\left[\lim_{r\to\infty}\frac 1r
\int_0^r(z_i(t+s)-\bar z^*_i)q_j(t)\dif t\right]\dif s,
\end{equation}
\begin{equation}
\BS R_{ijk}^{\hat H\to\bar z}=\int_0^\infty\left[\lim_{r\to\infty}\frac 1r
\int_0^r(z_i(t+s)-\bar z^*_i)q_j(t)q_k(t)\dif t\right]\dif s,
\end{equation}
\begin{equation}
\BS R_{ijk}^{\hat h\to\Sigma}=\int_0^\infty\left[\lim_{r\to\infty}
\frac 1r\int_0^r(z_i(t+s)-\bar z^*_i)(z_j(t+s)-\bar z^*_j)
q_k(t)\dif t\right]\dif s,
\end{equation}
\begin{equation}
\begin{split}
\BS R_{ijkl}^{\hat H\to\Sigma}=\int_0^\infty\bigg[\lim_{r\to\infty}
\frac 1r\int_0^r(z_i(t+s)-\bar z^*_i)(z_j(t+s)-\bar z^*_j)\times\\
\times q_k(t)q_l(t)\dif t-\Sigma^*_{ij}\delta_{kl}\bigg]\dif s.
\end{split}
\end{equation}
\end{subequations}
Returning back to the original perturbation coordinates, we obtain
\begin{subequations}
\label{eq:R_Gaussian_app}
\begin{equation}
\BS R_{ij}^{h\to\bar z}=\int_0^\infty\left[\lim_{r\to\infty}\frac 1r
\int_0^r(z_i(t+s)-\bar z^*_i)(z_k(t)-\bar z^*_k)\dif t\right]\dif s
\;\Sigma^{*-1}_{kj},
\end{equation}
\begin{equation}
\begin{split}
\BS R_{ijk}^{H\to\bar z}=\int_0^\infty\bigg[\lim_{r\to\infty}\frac 1r
\int_0^r(z_i(t+s)-\bar z^*_i)(z_l(t)-\bar z^*_l)\times\\\times
(z_k(t)-\bar z^*_k)\dif t\bigg]\dif s\;\Sigma^{*-1}_{lj},
\end{split}
\end{equation}
\begin{equation}
\begin{split}
\BS R_{ijk}^{h\to\Sigma}=\int_0^\infty\bigg[\lim_{r\to\infty}
\frac 1r\int_0^r(z_i(t+s)-\bar z^*_i)(z_j(t+s)-\bar z^*_j)\times\\
\times(z_l(t)-\bar z^*_l)\dif t\bigg]\dif s\;\Sigma^{*-1}_{lk},
\end{split}
\end{equation}
\begin{equation}
\begin{split}
\BS R_{ijkl}^{H\to\Sigma}=\int_0^\infty\bigg[\lim_{r\to\infty}
\frac 1r\int_0^r(z_i(t+s)-\bar z^*_i)(z_j(t+s)-\bar z^*_j)\times\\
\times(z_m(t)-\bar z^*_m)(z_l(t)-\bar z^*_l)\dif t\;
\Sigma^{*-1}_{mk}-\Sigma^*_{ij}\delta_{kl}\bigg]\dif s,
\end{split}
\end{equation}
\end{subequations}
which is given above in \eqref{eq:R_Gaussian}.\\

{\bf Acknowledgments.} The author is supported by the National Science
Foundation CAREER grant DMS-0845760, and the Office of Naval Research
grants N00014-09-0083 and 25-74200-F6607.

\end{document}